\documentclass[preprint,11pt]{article}
\usepackage{amssymb,amsfonts,amsmath,amsthm,amscd,dsfont,mathrsfs}
\usepackage{graphicx,float,psfrag,epsfig}
\usepackage{wrapfig}
\usepackage{algorithm}
\usepackage{algorithmic}
\usepackage{rotating}
\usepackage{multirow}

\DeclareMathAlphabet{\mathpzc}{OT1}{pzc}{m}{it}

\newtheorem{propo}{Proposition}[section]
\newtheorem{lemma}[propo]{Lemma}

\def\eps{\epsilon}

\def\reals{{\mathds R}}
\def\tM{\widetilde{M}}

\def\diag{{\rm diag}}

\def\E{\mathbb E}
\def\R{\mathbb R}

\def\ind{{\mathbb I}}

\def\<{\langle}
\def\>{\rangle}
\def\diag{{\rm diag}}
\def\Grass{{\sf G}}

\def\cP{{\cal P}}
\def\cF{{\cal F}}

\def\Xm{{\bf x}}
\def\Wm{{\bf w}}

\def\eps{\epsilon}

\def\hM{\widehat{M}}

\def\E{{\mathbb E}}
\def\grad{{\rm grad}\, }

\def\rank{{\rm rank}\, }
\def\hr{{\hat{r}}}
\def\Mmax{{M_{\rm max}}}
\def\optspace{{\sc OptSpace }}
\def\incoptspace{{\sc Incremental OptSpace }}

\begin{document}

\title{\optspace : A Gradient Descent Algorithm on the Grassman Manifold for Matrix Completion}
\author{Raghunandan~H.~Keshavan%,~\IEEEmembership{Student Member,~IEEE,}
        ~and~Sewoong~Oh
%\IEEEcompsocitemizethanks{
%\IEEEcompsocthanksitem R. H. Keshavan and S. Oh are with 
%the Departments of Electrical Engineering, Stanford University, Stanford, CA 94305.\protect\\
%E-mail: \{raghuram,swoh\}@stanford.edu% <-this % stops a space
%}
}

% make the title area
\maketitle
%
%============================================================
%
\begin{abstract}
We consider the problem of reconstructing a low rank 
matrix from a small subset of its entries. 
In this paper, 
we describe the implementation of
an efficient algorithm proposed in \cite{KOM09}, 
based on singular value decomposition followed by local manifold optimization, 
for solving the low-rank matrix completion problem.
It has been shown that if the number of revealed entries is large enough, 
the output of singular value decomposition gives a good estimate for the original matrix, 
so that local optimization reconstructs the correct matrix with high probability.
We present numerical results which show that 
this algorithm can reconstruct the low rank matrix exactly 
from a very small subset of its entries.
We further study the robustness of the algorithm
with respect to noise, and its performance on actual collaborative filtering datasets.
%Furthermore, we investigate the structure of real world data and 
%introduce novel modifications to our algorithm to account for the noise and human factors in the data.
%
\end{abstract}
%
%\begin{IEEEkeywords}
%Low-rank matrices, matrix completion, numerical analysis.
%\end{IEEEkeywords}
%\begin{center} \bfseries EDICS Category: MLR-LEAR \end{center}

%
%============================================================
%
\section{Introduction}
In this paper we consider the problem of reconstructing an $m\times n $ low rank matrix $M$ from a small set of observed
entries. This problem is of considerable practical interest and has many applications. One example
is collaborative filtering, where users submit rankings for small subsets of, say, movies, and the goal
is to infer the preference of unrated movies for a recommendation system \cite{Net06}. It is believed that
the movie-rating matrix is approximately low-rank, since only a few factors contribute to a user’s
preferences. Other examples of matrix completion include the problem of inferring 3-dimensional
structure from motion \cite{SFM} and triangulation from incomplete data of distances between wireless
sensors, also known as the sensor localization problem \cite{Singer08}, \cite{OKM09}.

%
%============================================================
%
\subsection{Prior and related work}\label{sec:relatedwork}

On the theoretical side, most recent work focuses on algorithms for exactly recovering the unknown
low-rank matrix. They provide an upper bound on the number of observed entries that guarantee
successful recovery with high probability. The main assumptions of this \emph{exact matrix completion}
problem are that the matrix $M$ to be recovered has rank $r \ll m,n$ and that the observed entries are
known exactly. The problem is equivalent to finding a minimum rank matrix matching the observed
entries . This problem is NP-hard. Adapting techniques from compressed sensing, Cand{\`e}s and Recht
introduced a convex relaxation of this problem \cite{CaR08}. They introduced the concept of incoherence
and proved that for a matrix $M$ of rank $r$ which has the incoherence property, solving the convex
relaxation correctly recovers the unknown matrix, with high probability, if the number of observed
entries $|E|$ satisfies, $|E| \geq C(\alpha)rn^{1.2}\log n$.

Recently \cite{KOM09} improved the bound to 
$|E| \geq C(\alpha)rn\max\{\log n, r\}$ for matrices with bounded condition number 
and provided an efficient algorithm called \optspace, based on spectral methods 
followed by local manifold optimization. 
For a bounded rank $r$, it is order optimal in the sense that
an $m \times n$ rank-$r$ matrix $M$ has $r(m+n-r)$ degrees of freedom
%the number of degrees of freedom in $M$ is $r(m+n-r)$ 
and without the same number of observations it is impossible to fix them. 
The extra $\log n$ factor is due to a coupon-collector effect:
it is necessary that $E$ contains at least one entry per row
and one per column, which happens only for $|E|\ge C n\log n$ \cite{CaR08,KMO08}.
 Cand{\`e}s and Tao proved a similar bound $|E| \geq C(\alpha)nr(\log n)^6$
with a stronger assumption on the original matrix $M$,
known as the \emph{ strong incoherence condition} \cite{CandesTaoMatrix} but without any assumption on the condition number of $M$.
For any value of $r$, it is only suboptimal by a poly-logarithmic factor\footnote{
In \cite{GLFBE09,Recht09,Gross09}, 
which appeared while we were preparing this manuscript, 
improved guarantees were proved for the convex relaxation algorithm. 
Namely, assuming only the incoherence property 
on the original matrix $M$,
the convex relaxation correctly recovers 
the matrix $M$ if $|E| \geq C(\alpha)nr(\log n)^2$.}. 
%Both results \cite{GLFBE09} 

When the pattern of observed entries is non-random, it is interesting to ask if the solution of the
rank $r$ matrix completion problem is unique \cite{Rigidity}. Building on the ideas from rigidity theory, 
Singer and Cucuringu introduce a randomized algorithm to determine whether it is possible to uniquely complete a
partially observed matrix to a matrix of specific rank $r$. Furthermore, by applying their algorithm to
random patterns of observed entries, one can get a lower bound on the minimum number of observed
entries necessary to correctly recover the matrix $M$.

While most theoretical work focuses on proving bounds for the exact matrix completion problem,
a more interesting and practical problem is when the matrix M is only approximately low rank or
when the observation is corrupted by noise. The main focus of this \emph{approximate matrix completion}
problem is to design an algorithm to find an $m \times n$ low-rank matrix $\hM$ that best approximates
the original matrix $M$ and provide a bound on the root mean squared error (RMSE) given by,
${\rm RMSE} = \frac{1}{\sqrt{mn}}||M-\hM||_F$. Cand{\`e}s and Plan introduced a generalization of 
the convex relaxation from \cite{CaR08} to the approximate case, 
and provided a bound on the RMSE \cite{CandesPlan}. 
More recently, a bound on the RMSE 
achieved by the \optspace algorithm with noisy observations was obtained in \cite{KMO09noise}. 
This bound is order optimal in a number of situations 
and improves over the analogous result in \cite{CandesPlan}.

 On the practical side, directly solving the convex relaxation introduced in \cite{CaR08} 
requires solving a Semidefinite Program (SDP), the complexity of which grows proportional to $n^3$. 
In the last year, many authors have proposed efficient algorithms 
for solving the low-rank matrix completion problem.
These include Singular Value Thresholding (SVT) \cite{CCS08},
Accelerated Proximal Gradient (APG) algorithm \cite{APG},
Fixed Point Continuation with Approximate SVD (FPCA) \cite{FPCA},
Atomic Decomposition for Minimum Rank Approximation ({\sc ADMiRA}) \cite{ADMiRA},
{\sc Soft-Impute} \cite{MHT09}, 
Subspace Evolution and
Transfer (SET) \cite{DM09}, 
Singular Value Projection (SVP) \cite{MJD09}
and {\sc OptSpace} \cite{KOM09}.

The problems each of these algorithms are trying to solve are described in Section \ref{sec:model}.
SVT is an iterative algorithm for 
solving the convex relaxation of the \emph{exact matrix completion} problem, 
which minimizes the nuclear norm (a sum of the singular values) 
under the constraints of matching the observed entries as in (\ref{P1}).
APG, FPCA and {\sc Soft-Impute} are efficient algorithms for solving the convex relaxation of the
\emph{approximate matrix completion} problem, which is a
nuclear norm regularized least squares problem as in (\ref{P3}).
{\sc ADMiRA} is an extension of Compressive Sampling Matching Pursuit (CoSaMP) \cite{CoSaMP}
, which is an iterative method for solving a least squares problem with bounded rank $r$ 
as described in (\ref{P4}). 
SVP is another approach to solve (\ref{P4}), which is a generalization of Iterative Hard Thresholding (IHT) \cite{BD09}.

%
%============================================================
%

\subsection{Contributions and outline}

The main contribution of this paper is to develop and implement efficient procedures, based on
the \optspace algorithm introduced in \cite{KOM09}, for solving the exact and approximate matrix
completion problems and add novel modifications, namely {\sc Rank Estimation} and {\sc Incremental
OptSpace}, that allow for a broader application and a better performance.

The algorithm described in \cite{KOM09} requires a knowledge of the rank of the original matrix. In
the following, we introduce a procedure, {\sc Rank Estimation}, that is guaranteed to correctly estimate
the rank of the original matrix from a partially revealed matrix under some conditions, which in turn
allows us to use the algorithm in a broader set of applications.

Next, we introduce {\sc Incremental OptSpace}, a novel modification to {\sc OptSpace}. We show
that, empirically, {\sc Incremental OptSpace} has substantially better performance than {\sc OptSpace}
when the underlying matrix is ill-conditioned.

%Further, we run numerical simulations under various scenarios and compare the performance of
%the algorithm, with respect to accuracy of reconstruction and speed, with the different algorithms.
Further, we carry out an extensive empirical comparison of various reconstruction algorithms.
This is particularly important because 
performance guarantees are only ``up to constants''
and therefore they have limited use in comparing different algorithms.
Finally, we apply our algorithm to real-world data and demonstrate that it is readily applicable to
the real data.

 The organization of the paper is as follows. In Section 2
we describe the low rank matrix completion problem and 
convex relaxations to the basic NP-hard approach, 
mostly to set our notation for later use.
Section 3 introduces an efficient implementation of the \optspace algorithm with novel modifications.
In Section 4 we discuss the results of numerical simulations with 
respect to speed and accuracy and compare the performance of \optspace with that of the other algorithms.

%The numerical simulations show that \optspace, 
%which has the best known performance bound 
%for both exact and approximate matrix completion problems, 
%also has the best performance in simulations with synthetic as well as real data.

%
%============================================================
%
\section{The model definition} \label{sec:model}

In the case of \emph{exact matrix completion}, 
we assume that the matrix $M$ has exact low rank $r\ll\min\{m,n\}$
and that the observed entries are known exactly. 
More precisely, we assume that 
there exist matrices $U$ of dimensions $m\times r$, 
V of dimensions $n \times r$, and a diagonal matrix $\Sigma$ 
of dimensions $r\times r$, such that 
\begin{eqnarray}
 M = U\Sigma V^T\; \label{eq:OriginalMatrix}.
\end{eqnarray}
Notice that for a given matrix $M$, the factors $(U,V,\Sigma)$ are not unique.

Out of the $m\times n$ entries of $M$, a subset $E\subseteq[m]\times[n]$ is observed.
Let $M^E$ be the $m \times n$ observed matrix that stores all the observed values, such that
\begin{eqnarray}
  M^E_{i,j} = \left\{
              \begin{array}{rl}
              M_{ij} & \text{if } (i,j)\in E\, ,\\
              0       & \text{otherwise.}
              \end{array} \right. \label{eq:ObservedMatrix}
\end{eqnarray}
Our goal is to find a low rank estimation $\hM(M^E,E)$ of the original matrix $M$ from
the observed matrix $M^E$ and the set of observed indices $E$.
 
If the number of observed entries $|E|$ is large enough, there is
a unique rank $r$ matrix which matches the observed entries.
In this case, solving the following optimization problem will recover the original matrix correctly.
\begin{eqnarray}
&\text{minimize}   & \rank(X)       \; \label{P0} \\
&\text{subject to} & \cP_E(X)=\cP_E(M)      \;, \nonumber
%&\text{subject to} & X_{ij}=M_{ij} \;,\;\, \forall(i,j)\in E \;, \nonumber
\end{eqnarray}
where $X\in\reals^{m\times n}$ is the variable matrix, $\rank(X)$ is the rank of matrix $X$, 
and $\cP_E(\cdot)$ is the projector operator defined as 
\begin{eqnarray}
\cP_E(M)_{ij} = \left\{\begin{array}{ll}
M_{ij} & \mbox{ if $(i,j)\in E$,}\\
0 & \mbox{otherwise.}
\end{array}\right.\label{eq:ProjectorDef}
\end{eqnarray}
The solution of this problem is the lowest rank matrix that matches the observed entries. Notice that
this is optimal in the sense that if this problem does not recover the correct matrix $M$ then there
%This problem finds the lowest rank matrix that matches the observed entries. 
%Notice that the solution of problem (\ref{P0}) is optimal. If this problem does not recover the correct matrix $M$ 
exists at least one other rank-$r$ matrix that matches all the observations, 
and no other algorithm can do better. However, this optimization problem is NP-hard and
all known algorithms require doubly exponential time in $n$ \cite{CaR08}.
This is especially inadequate since we are interested in cases where the dimension of the matrix $M$ 
is large ( eg. such as $5\cdot10^5 \times 2\cdot10^4$ for \cite{Net06}). 

In compressed sensing problems, minimizing the $\ell_1$ norm of a vector is used as a convex relaxation of 
minimizing the $\ell_0$ norm, or equivalently minimizing the number of non-zero entries, 
for sparse signal recovery. We can adopt this idea to matrix completion, where
$\rank(\cdot)$ of a matrix corresponds to $\ell_0$ norm of a vector, 
and nuclear norm to $\ell_1$ norm \cite{CaR08}, 
\begin{eqnarray}
&\text{minimize}   & ||X||_*       \; \label{P1}\\
&\text{subject to} & \cP_E(X)=\cP_E(M)      \;, \nonumber
\end{eqnarray}
where $||X||_*$ denotes the nuclear norm of $X$, i.e the sum of its singular values.

In the case of \emph{approximate matrix completion} problem, where
the observations are contaminated by noise or the original matrix 
to be reconstructed is only approximately low rank, the constraint $\cP_E(X)=\cP_E(M)$ 
must be relaxed. This results in either the problem \cite{CandesPlan,FPCA,APG,MHT09}
\begin{eqnarray}
&\text{minimize}   & ||X||_*       \; \label{P2}\\
&\text{subject to} & ||\cP_E(X)-\cP_E(M)||_F \leq \Theta      \;, \nonumber
\end{eqnarray}
or its Lagrangian version
\begin{eqnarray}
&\text{minimize}   & \mu||X||_*+\frac{1}{2}||\cP_E(X)-\cP_E(M)||_F^2    \;. \label{P3}
\end{eqnarray}
In \cite{ADMiRA,MJD09}, problem (\ref{P0}) is recast into the rank-$r$ matrix approximation problem of 
\begin{eqnarray}
&\text{minimize}   & ||\cP_E(X)-\cP_E(M)||_F  \; \label{P4}\\
&\text{subject to} & \rank(X)\leq r     \;. \nonumber
\end{eqnarray}

In the following, we present an efficient algorithm, namely {\sc OptSpace}, 
to solve the low-rank matrix completion problem which is closely related to (\ref{P4}),
and numerically compare its performance with those of the competing algorithms
in the case of exact as well as approximate matrix completion problems.

%
%======================================================
%
\section{Algorithm} \label{sec:Algorithm}

Algorithm 1 describes the overview of \optspace. 
Each step is explained in detail in the following sections.
The basic idea is to minimize 
the cost function $F:\reals^{m \times r}\times\reals^{n\times r}\to\reals$, defined as

\begin{eqnarray}
F(X,Y) 		&\equiv& \min_{S\in \reals^{r\times r}}\cF(X,Y,S)\, ,\label{eq:DefF1}\\
%\cF(X,Y,S)	&\equiv& \sum_{(i,j)\in E}f\Big(M_{ij},(XSY^T)_{ij}\Big) + \lambda_k\sum_{(i,j)\notin E}\frac{1}{2}\big(XSY^T\big)_{ij}^2 \, \label{eq:DefF2}.
\cF(X,Y,S)	&\equiv& \sum_{(i,j)\in E} f \Big(M_{ij},(XSY^T)_{ij}\Big) \, \label{eq:DefF2}.
\end{eqnarray}
Here $X\in \reals^{n\times r}$, $Y\in\reals^{m\times r}$ are 
orthogonal matrices, normalized as $X^{T}X = m\ind$, $Y^TY = n\ind$
where $\ind$ denotes the identity matrix, and
$f:\reals \times\reals\to\reals$ is a element-wise cost function. 
A common example that is useful in practice is the squared distance $f(x,y)=\frac{1}{2}(x-y)^2$.

Minimizing $F(X,Y)$ is an \emph{a priori} difficult task, since 
$F$ is a non-convex function. The basic idea is that the singular 
value decomposition (SVD) of $M^E$ provides an excellent initial guess,
and that the minimum can be found with high probability by
standard gradient descent after this initialization.
Two caveats must be added to this description: 
$(1)$ In general the matrix $M^E$ must be `trimmed' to eliminate
over-represented rows and columns; $(2)$ We need to estimate the target rank $r$.

\phantom{a} 

\vspace{0.3cm}
\begin{tabular}{ll} 
\hline
\vspace{-0.3cm}\\
\multicolumn{2}{l}{ {\bf Algorithm 1 :} {\sc OptSpace} }\\
\hline
\vspace{-0.3cm}\\
\multicolumn{2}{l}{{\bf Input:} observation matrix $M^E$, observed set $E$}\\
\multicolumn{2}{l}{{\bf Output:} estimated matrix $\hM$}\\
1: & Trim $M^E$, and let $\tM^E$ be the output;\\
2: & Estimate the rank of $M$, and let $\hr$ be the estimation;\\
3: & Compute the rank-$\hr$ projection of $\tM^E$, $\cP_\hr(\tM^E)=X_0S_0Y_0^T$;\\
4: & Minimize $F(X,Y)$ through gradient descent, with initial condition $(X_0,Y_0)$, \\
   & and return $\hM=XSY^T$.\\%, where $S$ is the minimizer defined in (Eq.~\ref{eq:DefF1})\\
\hline
\end{tabular}
\vspace{0.3cm}

%
%======================================================
%
\subsection{Trimming}

We say that a row is over-represented if its degree is more than $2|E|/m$ 
(twice the average degree), where degree of a row is defined as 
the number of observed entries in that row. Analogously, 
a column is over-represented if its degree is more than $2|E|/n$.
Trimming is a procedure that takes $M^E$ and $E$ as input and outputs $\tM^E$ 
by setting to $0$ all of the entries in over-represented rows and columns. 
Let $d_l(i)$ and $d_r(j)$ be the degree of $i^{th}$ row and $j^{th}$ column of $M$ respectively. 
Then the trimmed matrix $\tM^E$ is defined as 
\begin{eqnarray}
\tM^E_{ij} = \left\{
	\begin{array}{rl}
		0        & \text{if } d_l(i) > 2|E|/m \text{ or } d_r(j) > 2|E|/n\;,\\
    		M^E_{ij} & \text{ otherwise.}\\
	\end{array} \right.
\end{eqnarray}
The trimming step is essential when $|E| = \Theta(n)$, 
in which case there exists over-represented columns and rows of degrees $\Theta(\log n/\log\log n)$,
corresponding to singular values of the order $\Theta(\sqrt{\log n/\log\log n})$. 
As $n$ grows large, while these spurious singular values dominate 
the principal components in step 3 of the Algorithm 1, 
the corresponding singular vectors are highly concentrated on the 
over-represented rows and columns (respectively for left and right singular vectors) 
and do not provide any useful information about the unobserved entries of $M$.
\subsection{Estimating the rank}
\label{sec:rankestimation}
Define $\eps \equiv |E|/\sqrt{mn}$. 
In the case of a square matrix $M$, 
$\eps$ corresponds to the average degree per row or per column. 
Throughout this paper, the parameter $\eps$ will be frequently
used as the model parameter indicating how difficult the problem instance is.

By singular value decomposition of the trimmed matrix, let
\begin{equation}
\tM^E = \sum_{i=1}^{\min(m,n)}\sigma_i  x_iy_i^T\, , \label{eq:MatrixSVD}
\end{equation}
where $x_i$ and $y_i$ are the left and right singular vectors 
corresponding to $i$th singular value $\sigma_i$.
Then, the following cost function is defined in terms of the singular values.
\begin{eqnarray}
 R(i)=\frac{\sigma_{i+1}+\sigma_1\sqrt{\frac{i}{\eps}}}{\sigma_i} \;.
\end{eqnarray}

%\begin{algorithm}
% \caption{Rank Estimation}\label{alg:rankestimation}
% \begin{algorithmic}[1]
%  \REQUIRE trimmed observation matrix $\tM^E$
%  \ENSURE estimated rank $\hr$
%  \STATE SVD $\tM^E$, and let $\{\sigma_i\}$ be the set of singular values;
%  \STATE Find rank $r$ that minimizes $R(i)=\frac{\sigma_{i+1}+\sqrt{i\eps}}{\sigma_i}$.
% \end{algorithmic}
%\end{algorithm}
Based on the above definition, {\sc Rank Estimation} consists of two steps:

\vspace{0.3cm}
\begin{tabular}{ll} 
\hline
\multicolumn{2}{l}{ {\sc Rank Estimation}}\\
\hline
\vspace{-0.3cm}\\
\multicolumn{2}{l}{{\bf Input:} trimmed observation matrix $\tM^E$}\\
\multicolumn{2}{l}{{\bf Output:} estimated rank $\hr$}\\
%\hline
%\vspace{-0.4cm}\\
1: & Compute singular values $\{\sigma_i\}$ of $\tM^E$;\\
2: & Find the index $i$ that minimizes $R(i)$, and let $\hr$ be the minimizer.\\
\vspace{-0.3cm}\\
\hline
%\label{alg:rank}
\end{tabular}
\vspace{0.2cm}

The idea behind this algorithm is that, if enough entries of $M$ are revealed 
then there is a clear separation between the first $r$ singular values, 
which reveal the structure of the matrix $M$ to be reconstructed, 
and the spurious ones \cite{KOM09}. 
As described in the following proposition, 
%assuming that each entries of $M$ is observed independently with probability $\eps/\sqrt{mn}$,
we can show that this simple procedure is guaranteed to 
reconstruct the correct rank $r$, with high probability, for $|E|$ large enough. 
For the proof of this proposition, we refer to Appendix \ref{app:rank}.
\begin{propo}\label{pro:rank}
Assume $M$ to be a rank $r$ $m\times n$ matrix with 
bounded condition number $\kappa$.
Then there exists a constant $C(\kappa)$ such that,
if $\eps>C(\kappa)r$, then {\sc Rank Estimation} correctly estimates the rank $r$, with high probability.
\end{propo}
%

%
%======================================================
%
\subsection{Rank-$\rho$ projection}

Rank-$\rho$ projection consists of performing a sparse SVD on 
$\tM^E$ and rescaling the singular values and singular vectors appropriately.
From the {\sc Rank Estimation} step we have the SVD of $\tM^E$ in Eq.~(\ref{eq:MatrixSVD}), 
namely $\tM^E = \sum_{i=1}^{\min(m,n)}\sigma_i  x_iy_i^T$. 
Define the projection : 
\begin{eqnarray}
 \cP_\rho(\tM^E)=X_0S_0Y_0^T \;, \label{eq:projection}
\end{eqnarray}
for normalized orthogonal matrices $X_0\in\reals^{m \times \rho}$ and 
$Y_0\in\reals^{n\times \rho}$, and a $\rho \times \rho$ diagonal matrix $S_0$,
defined in terms of the singular values and singular vectors in Eq.~(\ref{eq:MatrixSVD})
as $X_0 = \sqrt{m}[x_1,\ldots,x_\rho]$, $Y_0 = \sqrt{n}[y_1,\ldots,y_\rho]$,
and $S_0=({1}/{\eps})\text{diag}(\sigma_1,\ldots,\sigma_\rho)$. Notice that we 
do not need to compute the scaled singular values $S_0$, since we only 
require $X_0$ and $Y_0$ for the following local optimization step.
There are a number of low complexity algorithms available for forming a sparse SVD,
as well as a number of open source implementations of these algorithms.

%\vspace{0.3cm}
%
%\begin{tabular}{ll} 
%\hline
%\multicolumn{2}{l}{ {\sc Projection} ( matrix $\tM^E$, rank $r$ )}\\
%\hline
%1: & Compute $\tX$ and $\tY$, the left and right singular matrices of $\tM^E$ upto rank $r$ \\
%2: & return $\Xm_0 = (\sqrt{m} \tX, \sqrt{n} \tY)$ \\
%\hline
%\end{tabular}
%
%\vspace{0.3cm}

%
%======================================================
%
\subsection{Gradient descent on the Grassman manifold}

The {\sc Manifold Optimization} step involves gradient descent with 
variables $X\in\reals^{m\times r}$ and $Y\in\reals^{n\times r}$ 
using the cost function $F(X,Y)$ defined below. 
In this section, we use $r$ and $\hr$ interchangeably to denote the estimated rank of matrix $M$.
\begin{eqnarray}
F(X,Y) 		&\equiv& \min_{S\in \reals^{r\times r}}\cF(X,Y,S)\, ,\label{eq:DefF5}\\
\cF(X,Y,S)	&\equiv& \sum_{(i,j)\in E}f\Big(M_{ij},(XSY^T)_{ij}\Big) + \lambda \sum_{(i,j)\notin E}\frac{1}{2}\big(XSY^T\big)_{ij}^2 \, \label{eq:DefF3},
\end{eqnarray}
where $f:\reals \times\reals\to\reals$ is an element-wise cost function. 
Note that compared to Eq.~(\ref{eq:DefF2}), we have additional term in Eq.~(\ref{eq:DefF3}),
which is a regularization term with a regularization coefficient $\lambda \in [0,1]$.

 The above general formulation allows for a freedom in choosing a suitable 
cost function $f$ for different applications.
% such as the discrete movie ranking matrix of the Netflix Data.
However, a common example of the cost function $f(x,y)=\frac{1}{2}(x-y)^2$
works very well in practice as well as in proving performance bounds \cite{KOM09}.
Hence, throughout this paper, we use the squared difference as the cost function, 
resulting in 
\begin{eqnarray*}
\cF(X,Y,S) \equiv \frac{1}{2}\Big|\Big|\cP_E(M-XSY^T)\Big|\Big|_F^2 + \lambda \frac{1}{2}\Big|\Big|\cP_{E^\perp}(XSY^T)\Big|\Big|_F^2 \;, %\label{eq:DefF4}
\end{eqnarray*}
where the projector operator $\cP_E$ for a given $E$ is defined in Eq.~(\ref{eq:ProjectorDef}), 
and $E^\perp$ is the complementary set of $E$.

For the results in this paper, we choose $\lambda = 0$ but we observe that using a positive $\lambda$ helps
when the matrix entries are corrupted by noise. For $\lambda = 0$, the gradient of $F(X, Y )$ can be written
explicitly as

\begin{eqnarray*}
\grad F(\Xm)_{X} & = & \mathcal{P}_E(XSY^T - M)YS^T \;,\\
\grad F(\Xm)_{Y} & = & \mathcal{P}_E(XSY^T - M)^TXS \;,
\end{eqnarray*}
where $S$ is the $r \times r$ matrix that achieves the minimum in the definition of $F(X,Y)$, Eq.~(\ref{eq:DefF5}).

One important feature of \optspace is that 
$F(X,Y)$ is regarded as a function
of the $r$-dimensional subspaces of $\reals^m$ and $\reals^n$
generated (respectively) by the columns of $X$ and $Y$.
This interpretation is justified by the fact that
$F(X,Y) = F(XA,YB)$ for any two orthogonal matrices 
$A$, $B\in \reals^{r\times r}$.
The set of $r$ dimensional subspaces of $\reals^m$ is 
a differentiable Riemannian manifold $\Grass(m,r)$ (the Grassman
manifold) \cite{KOM09}. 
The gradient descent algorithm is applied to the function $ F: \Grass(m,r)\times\Grass(n,r) \to  \reals \;$ 
For further details on optimization by gradient descent on matrix manifolds we
refer to \cite{Edelman,ManifBook}.

In the following, we use a compact representation $\Xm $ for a pair $(X,Y)$, 
with $X$ an $n\times r$ matrix and $Y$ an $m\times r$ matrix. Similarly, 
the gradient is represented by : $\grad F(\Xm_k)=(\grad F(\Xm_k)_X,\grad F(\Xm_k)_Y)$.
Let $\Xm_0 = (X_0,Y_0)$, where $X_0$ and $Y_0$ are the 
normalized left and right singular matrices from rank-$r$ projection. 
The {\sc Manifold Optimization} algorithm starting at $\Xm_0$ is described below. 
We refer to \cite{KOM09} for justifications and performance bounds of the algorithm.

For any scalar $\tau$, it is shown in \cite{Armijo} that 
this algorithm converges to the local minimum. 
However, numerical experiments suggest $\tau=10^{-3}$ is a good choice. 
The algorithm stops when the fit error $||\cP_E(M-\hM)||_F/||\cP_E(M)||_F$
goes below some threshold $\delta_{\rm tol}$, e.g. $10^{-6}$. 
The basic idea is that this is a good indicator of the relative error on the whole set, 
$||M-\hM||_F/||M||_F$. This stopping criterion is also used in 
other algorithms such as SVT in \cite{CCS08} where the authors provide a convincing argument for its use.

\vspace{.3cm}
\begin{tabular}{rl}
\hline
\vspace{-.3cm}\\
\multicolumn{2}{l}{ {\sc Manifold Optimization}}\\
\hline
\vspace{-.3cm}\\
\multicolumn{2}{l}{ {\bf Input:} observed matrix $M^E$, estimated rank $\hr$, initial factors $\Xm_0=(X_0,Y_0)$,} \\
\multicolumn{2}{l}{ $\;\;\;\;\;\;\;\;\;\;\;\;$ tolerance $\delta_{\rm tol}$, maximum iteration count $k_{max}$, step size $\tau$ }\\
\multicolumn{2}{l}{ {\bf Output:} reconstructed matrix $\hM$ }\\
\vspace{-.4cm}\\
1: & {\bf For} $k=0,1,\dots,k_{max}$ {\bf do}: \\
%2: &\hspace{0.2cm} Compute $\Wm_k = \grad F(\Xm_k) = (\grad F(\Xm_k)_X,\grad F(\Xm_k)_Y)$;\\
2: &\hspace{0.2cm} Compute $S_k=\arg \min_S\{\mathcal{F} (X_k,Y_k,S)\}$\\
3: &\hspace{0.2cm} Compute $\Wm_k = \grad F(\Xm_k)$\\
%4: &\hspace{0.2cm} Let $\Xm_k(t) = \Xm_k+\Wm_k t$\\
%5: &\hspace{0.2cm} Find $t_{k}^{*}$ using backtracking line search\\
4: &\hspace{0.2cm} Set $t_{k}=\tau$\\
5: &\hspace{0.2cm} {\bf While} $F(\Xm_k- t_{k} \Wm_k) - F(\Xm_k) > \frac{1}{2} t_{k} ||\Wm_k||^2$, {\bf do}\\
6: &\hspace{0.6cm} $t_{k} \leftarrow t_{k}/2$ \\
7: &\hspace{0.2cm} Set $\Xm_{k+1} = \Xm_k- t_k\Wm_k$ \\
8: &\hspace{0.2cm} {\bf If } $||\cP_E(M-\hM)||_F/||\cP_E(M)||_F  < \delta_{\rm tol}$ {\bf then break}\\
9: & {\bf End for}\\
10:& Set $\hM=X_kS_kY_k^T$\\
\hline
\end{tabular}

\vspace{.3cm}

%We now describe the backtracking algorithm used in Step 5. of {\sc Gradient Descent}. 
%\vspace{.3cm}
%\begin{tabular}{ll}
%\hline
%\multicolumn{2}{l}{ {\sc Backtracking line search} ( matrix $M^E$, factors
%$\Xm$, gradient $\Wm$, scalar $\alpha$ )}\\
%\hline
%1: &  While $F(\Xm- \alpha \Wm) - F(x) > \frac{1}{2} \alpha ||\Wm||^2$, do \\
%2: &\hspace{0.2cm} $\alpha \leftarrow \alpha/2$ ;\\
%3: & return $\alpha$;\\
%\hline
%\end{tabular}
%\vspace{.3cm}

%
%======================================================
%
\subsection{A novel modification to \optspace for ill-conditioned matrices}
\label{sec:IncrementalOptSpace}

In this section, we describe a novel modification to the \optspace algorithm, 
%with added complexity, 
which has substantially better performance 
in the case when the matrix $M$ to be reconstructed is ill-conditioned. 
%the modified \optspace algorithm for matrices with high condition number. 
When the condition number $\kappa(M)$ is high, 
the initial guess in step 3 of \optspace 
for ($u_r,v_r$), the singular vectors which correspond to the smallest singular value,  
are often far from the correct ones.
To compensate for this discrepancy, 
we start by first finding ($u_1,v_1$), 
the singular vectors corresponding to the first singular value,
and incrementally search for the next ones.
%The general idea consists of optimizing the space incrementally.

\vspace{0.3cm}
\begin{tabular}{ll} 
\hline
\vspace{-0.3cm}\\
\multicolumn{2}{l}{ {\bf Algorithm 2 :} {\sc Incremental OptSpace} }\\
\hline
\vspace{-0.3cm}\\
\multicolumn{2}{l}{{\bf Input:} observation matrix $M^E$, observed set $E$, }\\
\multicolumn{2}{l}{ $\;\;\;\;\;\;\;\;\;\;\;\;$ tolerance $\delta_{\rm tol}$, maximum rank count $\rho_{max}$}\\
\multicolumn{2}{l}{{\bf Output:} estimation $\hM$}\\
1: & Trim $M^E$, and let $\tM^E$ be the output\\
2: & Set $\hM^{(0)} = 0$\\
3: & {\bf For} $\rho=0,1,\dots,\rho_{\rm max}$ {\bf do}: \\
4: & \hspace{0.2cm}Compute the rank-$1$ projection of $\tM^E - \hM^{(\rho)}$, $\cP_\hr(\tM^E - \hM^{(\rho)}) = X_0^{(\rho)}S^{(\rho)}_0 Y^{(\rho)T}_0$\\
5: & \hspace{0.2cm}Set $X^{(\rho)}_0 =  [X^{(\rho-1)} ;\, X^{(\rho)}_0; ]$ and 
     $Y^{(\rho)}_0 =  [Y^{(\rho-1)} ;\, Y^{(\rho)}_0; ]$ \\
6: & \hspace{0.2cm}Minimize $F(X,Y)$ through {\sc Manifold Optimization} with $\rho$ replacing $\hr$,\\ 
   & \hspace{0.2cm} with initial condition $(X^{(\rho)}_0,Y^{(\rho)}_0)$\\
   & \hspace{0.2cm} and stopping criterion of $|F(\Xm_{k+1}) - F(\Xm_k)| \le \delta_{\rm tol}F(\Xm_k)$, \\
  &  \hspace{0.2cm} and let $\hM^{(\rho)} = X^{(\rho)}S^{(\rho)}Y^{(\rho)T}$ be the output\\
7: & \hspace{0.2cm}{\bf If } $||\cP_E(M-\hM^{(\rho)})||_F/||\cP_E(M)||_F  < \delta_{\rm tol}$ {\bf then break}\\
9: & {\bf End for}\\
10: & Return $\hM^{(\rho)}$.\\
\hline
\end{tabular}
\vspace{0.3cm}
% 
% \begin{tabular}{rl}
% \hline
% \vspace{-.4cm}\\
% \multicolumn{2}{l}{ {\sc Modified Manifold Optimization}}\\
% \hline
% \vspace{-.4cm}\\
% \multicolumn{2}{l}{ {\bf Input:} observed matrix $M^E$, rank r, initial factors $\Xm_0=(X_0,Y_0)$,} \\
% \multicolumn{2}{l}{ $\;\;\;\;\;\;\;\;\;\;\;\;$ tolerance $\delta_{tol}$, maximum iteration count $k_{max}$, step size $\tau$ }\\
% \multicolumn{2}{l}{ {\bf Output:} reconstructed matrix $\hM = XSY^T$ }\\
% \vspace{-.4cm}\\
% 1: & {\bf For} $k=0,1,\dots,k_{max}$ {\bf do}: \\
% %2: &\hspace{0.2cm} Compute $\Wm_k = \grad F(\Xm_k) = (\grad F(\Xm_k)_X,\grad F(\Xm_k)_Y)$;\\
% 2: &\hspace{0.2cm} Compute $S_k=\arg \min_S\{\mathcal{F} (X_k,Y_k,S)\}$\\
% 3: &\hspace{0.2cm} Compute $\Wm_k = \grad F(\Xm_k)$\\
% %4: &\hspace{0.2cm} Let $\Xm_k(t) = \Xm_k+\Wm_k t$\\
% %5: &\hspace{0.2cm} Find $t_{k}^{*}$ using backtracking line search\\
% 4: &\hspace{0.2cm} Set $t_{k}=\tau$\\
% 5: &\hspace{0.2cm} {\bf While} $F(\Xm_k- t_{k} \Wm_k) - F(\Xm_k) > \frac{1}{2} t_{k} ||\Wm_k||^2$, {\bf do}\\
% 6: &\hspace{0.6cm} $t_{k} \leftarrow t_{k}/2$ \\
% 7: &\hspace{0.2cm} Set $\Xm_{k+1} = \Xm_k- t_k\Wm_k$ \\
% 8: &\hspace{0.2cm} {\bf If } $|F(\Xm_{k+1}) - F(\Xm_{k})|  < \delta_{tol} F(\Xm_k)$ {\bf then break}\\
% 9: & {\bf End for}\\
% 10:& Set $\hM=X_kS_kY_k^T$\\
% \hline
% \end{tabular}
% 
% \vspace{.3cm}

In the following numerical simulations, we demonstrate that 
{\sc Incremental OptSpace} brings significant performance gains 
when applied to ill-conditioned matrices, cf. Section \ref{sec:synthetic}.

%
%==============================================================
%
%\section{Numerical results}\label{sec:Implementation}
%
%======================================================
%
%\subsection{Implementation details} \label{sec:detail}

%
%======================================================
%
\section{Numerical results with randomly generated matrices}
\label{sec:synthetic}

 The \optspace algorithm described above 
was implemented in C\footnote{The code is available at http://www.stanford.edu/$\sim$raghuram/optspace.html} and tested on a 3.4 GHz Desktop computer with 4 GB RAM. 
For efficient singular value decomposition of sparse matrices, 
we used (a modification of) SVDLIBC\footnote{Available at http://tedlab.mit.edu/$\sim$dr/svdlibc/} 
which is based on SVDPACKC. 
In this section, we compare the performance of \optspace 
with other algorithms by numerical simulations.
In Section \ref{sec:exactcompletion}, 
the algorithms are tested on randomly generated matrices 
with noiseless observations, and in Section \ref{sec:noisycompletion} we compare the algorithms 
when we have noisy observations under different scenarios. 
%and in Section \ref{sec:real}, we compare the performance when applied to real data. 
% \textbf{[Apparently we need to send them a copy of the paper]}

For exact matrix completion experiments, 
we use $n\times n$ test matrices $M$ of rank $r$ 
generated as $M = UV^T$. Here, $U$ and $V$ are $n\times r$ matrices 
with each entry being sampled independently from a standard 
Gaussian distribution ${\cal N}(0,1)$, unless specified otherwise.
Then, each entry is revealed independently with probability $\eps/n$, so that 
on an average $n\eps$ entries are revealed. 
Numerical results show that there is no notable difference 
if we choose the revealed set of entries $E$ 
uniformly at random over all the subsets of the same 
size $|E|=n\eps$. 
We use $\delta_{\rm tol} = 10^{-5}$ and  $k_{max} = 1000$ as the stopping criteria. 

For approximate matrix completion, the matrices are generated as above and 
corrupted by additive noise $Z_{ij}$. 
First, in the standard scenario, $Z_{ij}$'s are 
independently and identically distributed 
according to a Gaussian distribution. 
For comparison, we also present 
numerical simulation results with different types of noise 
in the following subsections.
Again, each entry is revealed independently 
with a probability $\eps/n$. We use 
$||\cP_E(\hM - (M+Z) )||_F^2 \le (1 + \eps)|E| \sigma_n^2$ \cite{CCS08} 
(where $\sigma_n^2$ is the noise variance per entry) as the stopping criterion.

%
%======================================================
%
\subsection{Exact matrix completion}
\label{sec:exactcompletion}
\begin{figure}
\begin{center}
%\hspace{-2.1cm}
\includegraphics[width=10cm]{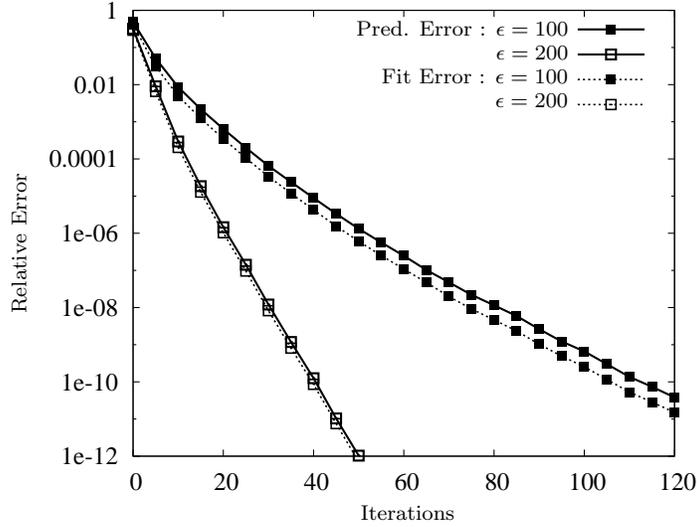}
\put(-125,160){\scriptsize{Fit Error : $\epsilon = 100$}}
\put(-125,151){\scriptsize{\phantom{Fit Error : }$\epsilon = 200$}}
\put(-134.5,178){\scriptsize{Pred. Error : $\epsilon = 100$}}
\put(-134.5,169){\scriptsize{\phantom{Pred. Error : }$\epsilon = 200$}}
\put(-133,-5){\scriptsize{Iterations}}
\put(-265,73){\begin{sideways}\scriptsize{Relative Error}\end{sideways}}
\end{center}
\caption{{\small Prediction and fit errors versus the number of iterations of the {\sc Manifold Optimization} step for rank $10$ matrices of dimension $n\times n$ with $n=1000$. Each entry is sampled with probability $\eps/n$ for two different values of $\eps$: $100$ and $200$.}} \label{fig:convergence}
%\vspace{-1.cm}
\end{figure}
We first illustrate the rate of convergence of \optspace. 
In Figure \ref{fig:convergence}, we plot the fit error, 
$||\cP_E(\hM-M)||_F/n$ and the prediction error 
$||\hM-M||_F/n$, with respect to the number of 
iterations of the {\sc Manifold Optimization} step. 
These plots are obtained for matrices with $n = 1000$ 
and $r = 10$ and averaged over $10$ instances. 
The results are shown for two values of $\eps$: $100$ and $200$.
We can see that the prediction error decays 
exponentially with the number of iterations in both cases.
Also, the prediction error is very close to the fit error, thus 
lending support to the validity of the chosen stopping criterion.
%Also, the prediction error is within a constant factor of the fit error, 
%thus lending support to the validity of the chosen stopping criterion.

We next study the {\em reconstruction rate} of the algorithm. We declare a matrix to be reconstructed if $||M - \hM||_F / ||M||_F \leq 10^{-4}$. The reconstruction rate is the fraction of instances for which the matrix was reconstructed.

\begin{figure}
%\hspace{-1.5cm}
\begin{center}
\includegraphics[width=10cm]{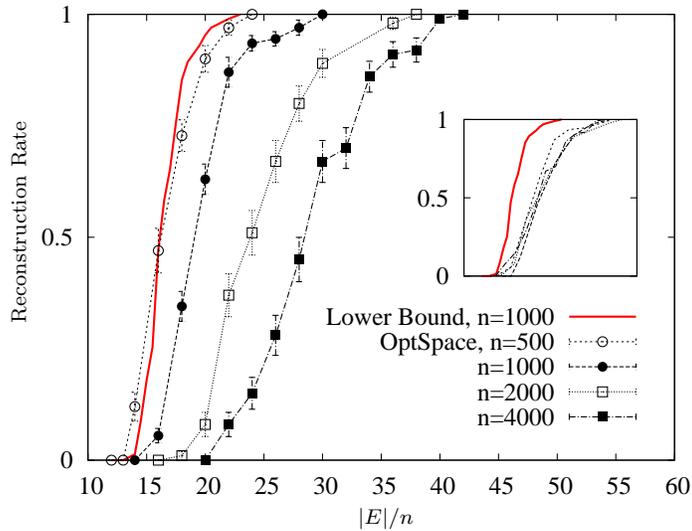}
\put(-134,-5){\scriptsize{$|E|/n$}}
\put(-265,71){\begin{sideways}\scriptsize{Reconstruction Rate}\end{sideways}}
\caption{{\small Reconstruction rates for rank $4$ matrices using \optspace for different matrix sizes. The solid curve is bound proved in \cite{Rigidity}. In the inset, the same data are plotted vs. $|E|/(n(\log n)^2)$ }} \label{fig:recon_n}
\end{center}
\end{figure}

In Figure \ref{fig:recon_n}, we plot the reconstruction rate as 
function of $|E|/n$ for \optspace on randomly generated rank-$4$ 
matrices for different matrix sizes $n$. 
As predicted by Theorem 1.2 of \cite{KOM09}, 
threshold of the reconstruction rate of {\sc OptSpace} 
is upper bounded by $|E|=Cn (\log n)^2$, for fixed rank $r=4$.
Here, an extra factor of $\log n$ comes from the fact that 
if we generate random factors $U$ and $V$ from a Gaussian distribution, 
then the incoherence parameter $\mu_0$ scales like $\log n$. 
However, the location of the threshold is surprisingly close to 
the lower bound proved in \cite{Rigidity} which scales as $|E|=Cn \log n$.
The lower bound provides a threshold below which the problem admits more than one solution.
Note that the lower bound is displayed only for the case when $n=1000$. 

\begin{figure}
\begin{center}
\includegraphics[width=10cm]{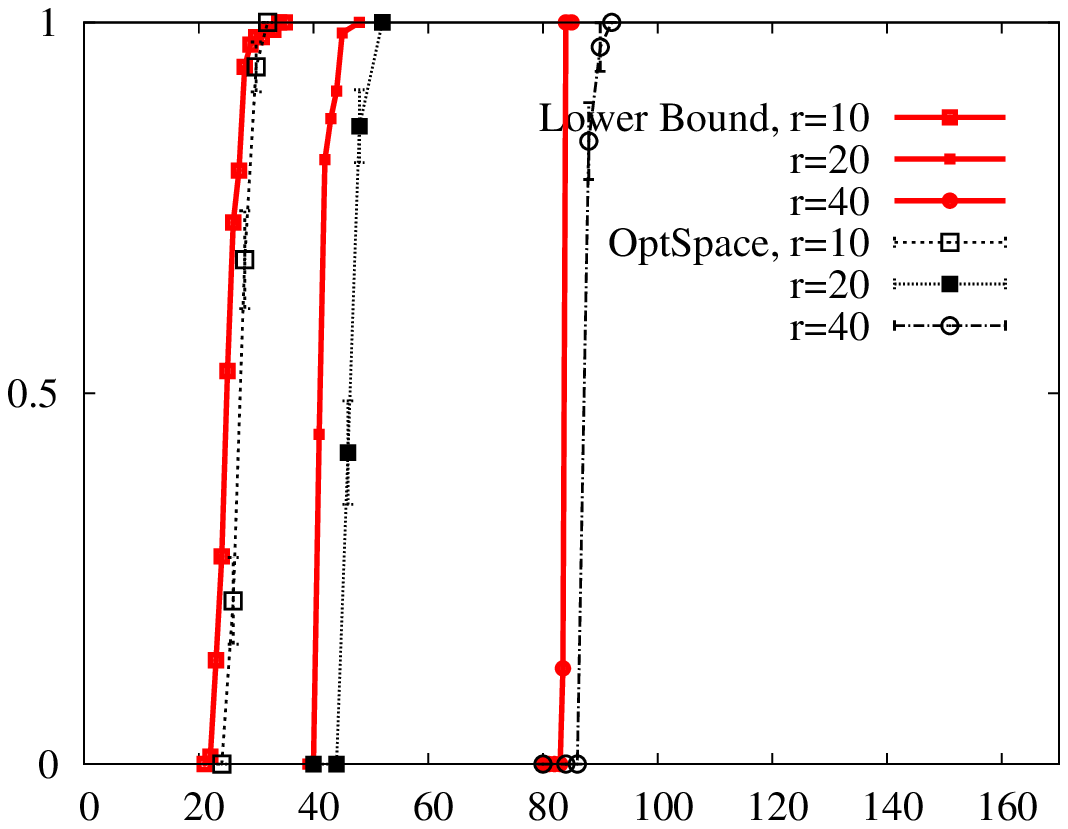}
%\hspace{-1.5cm}
%\input{n500.tex}
\put(-134,-5){\scriptsize{$|E|/n$}}
\put(-265,71){\begin{sideways}\scriptsize{Reconstruction Rate}\end{sideways}}
\caption{{\small Reconstruction rates for matrices with dimension $m=n=500$ using \optspace for different ranks. The solid curves are the bounds proved in \cite{Rigidity}. }} \label{fig:recon_r}
\end{center}
\end{figure}

In Figure \ref{fig:recon_r}, we plot the reconstruction rate for 
randomly generated matrices with dimensions $m=n=500$ 
using \optspace.
The resulting reconstruction rate is plotted for different ranks $r$ 
as a function of $|E|/n$.
As rank increases and for fixed $n$, 
the reconstruction rate has a sharp threshold at $|E|=Crn\log n$.
This indicates that in practice the dependence of the threshold 
on the rank scales like $r$ rather than $r^2$ as predicted by 
Theorem 1.2 of \cite{KOM09}.
Also, for all values of rank, the location of the threshold is 
surprisingly close to the lower bound proved in \cite{Rigidity}, 
below which the problem admits more than one solution. 

\begin{figure}
\begin{center}
\includegraphics[width=10cm]{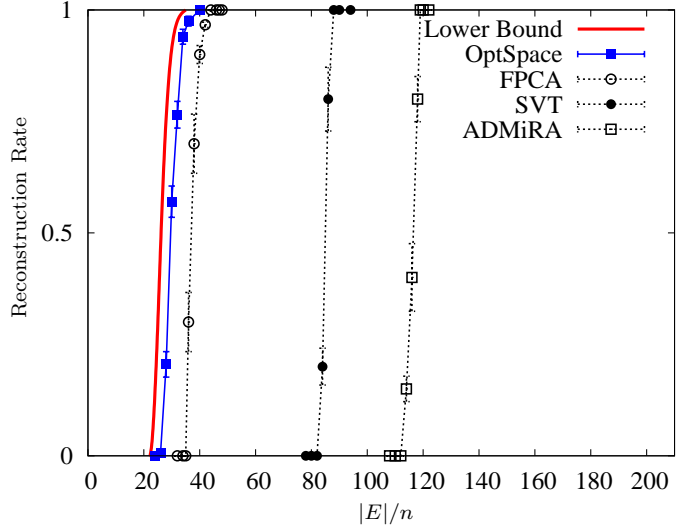}
%\hspace{-1.5cm}
%\input{comparisonexact.tex}
\put(-134,-5){\scriptsize{$|E|/n$}}
\put(-265,71){\begin{sideways}\scriptsize{Reconstruction Rate}\end{sideways}}
\caption{{\small Comparison of reconstruction rates for randomly generated rank $10$ matrix of dimension $m=n=1000$ for the \optspace and competing algorithms : {\sc FPCA}, {\sc SVT} and {\sc ADMiRA} \cite{CCS08,FPCA,ADMiRA}. The leftmost solid curve is a lower bound proved in \cite{Rigidity}.}} \label{fig:recon_rate}
\vspace{-1.cm}
\end{center}
\end{figure}

In Figure \ref{fig:recon_rate}, we plot the reconstruction rate of \optspace as a function of 
$|E|/n$ for rank $10$ matrices of dimension $m=n=1000$. 
Also plotted are the reconstruction rates obtained for the convex relaxation approach of 
\cite{CandesTaoMatrix} solved using the Singular Value Thresholding algorithm \cite{CCS08}, 
the FPCA algorithm from \cite{FPCA} and {\sc ADMiRA} \cite{ADMiRA}. 
We compare these with a theoretical lower bound on the reconstruction rate described in \cite{Rigidity}. 
Various algorithms exhibit threshold at different values of $|E|/n$, 
and the threshold depends on the problem size $n$ and the rank $r$.
This figure clearly illustrates that \optspace outperforms the other algorithms
on random data, and this was consistent for various values of $n$ and $r$.

In the following Tables \ref{tab:hard} and \ref{tab:easy}, 
we present numerical results obtained using these algorithms for 
different values of $n$ and $r$. 
Table \ref{tab:hard} presents results for smaller values of $\eps$ and 
hence for {\em hard} problems, whereas Table \ref{tab:easy} presents 
results for larger values of $\eps$ which are relatively {\em easy} problems. 
Note that the values of $\eps$ used in Table \ref{tab:hard} all correspond to 
$|E| \le 2.6d(n,r)$ where $d(n,r) = 2nr - r^2$ is the number of degrees of freedom. We ran into \emph{Out of Memory} problems for the FPCA algorithm for $n \ge 20000$ and hence we omit these problems from the table. All the results presented in tables are averaged over 5 instances.
On the easy problems, all the algorithms achieved similar performances, whereas on the hard problems, 
\optspace outperforms other algorithms on most of instances.

{\renewcommand{\arraystretch}{1.5}
\renewcommand{\tabcolsep}{0.2cm}}

%\begin{figure*}
\begin{table*}
\begin{center}
\footnotesize

\begin{tabular}{c c c | c c | c c | c c | c c }
\hline
\multirow{2}{*}{$n$} & \multirow{2}{*}{$r$} & \multirow{2}{*}{$\epsilon$} & \multicolumn{2}{|c|}{\optspace} & \multicolumn{2}{|c|}{SVT}& \multicolumn{2}{|c|}{FPCA}& \multicolumn{2}{|c}{{\sc ADMiRA}}\\
& & & rel. error & time(s) & rel. err. & time(s) & rel. err. & time(s) & rel. err. & time(s) \\
\hline
\multirow{3}{*}{$1000$} & $10$ & $50$ &
$1.95 \times 10^{-5}$ & $33$ &
$3.42 \times 10^{-1}$ & $734$ &
$6.04 \times 10^{-4}$ & $65$ & 
$4.41 \times 10^{-1}$ & $8 $\\

& $50$ & $200$ & 
$1.28 \times 10^{-5}$ & $235$ &
$2.54 \times 10^{-1}$ & $11769$ &
$1.07 \times 10^{-5}$ & $83$  &
$3.54 \times 10^{-1}$ & $141 $\\

& $100$ & $400$ &
$9.22 \times 10^{-6}$ & $837$ &
$7.99 \times 10^{-2}$ & $27276$ &
$3.86 \times 10^{-6}$ & $165$  &
$1.28 \times 10^{-1}$ & $1767 $\\

\hline

\multirow{3}{*}{$5000$} & $10$ & $50$ &
$7.27 \times 10^{-5}$ & $338$ &
$5.34 \times 10^{-1}$ & $476$ &
$9.99 \times 10^{-1}$ & $1776$ &
$5.13 \times 10^{-1}$ & $77 $\\

& $50$ & $200$ &
$1.47 \times 10^{-5}$ & $1930$ &
$4.87 \times 10^{-1}$ & $36022$ & 
$2.17 \times 10^{-2}$ & $2757$ &
$5.36 \times 10^{-1}$ & $358 $\\

& $100$ & $400$ &
$1.38 \times 10^{-5}$ & $6794$ &
$4.12 \times 10^{-1}$ & $249330$ &
$2.49 \times 10^{-5}$ & $3942$ &
$4.84 \times 10^{-1}$ & $36266 $\\

\hline
\multirow{3}{*}{$10000$} & $10$ & $50$ &
$1.91 \times 10^{-5}$ & $725$ &
$6.33 \times 10^{-1}$ & $647$ &
$9.99 \times 10^{-1}$ & $9947$ &
$6.19 \times 10^{-1}$ & $129 $\\

& $50$ & $200$ & 
$5.02 \times 10^{-6}$ & $3032$ &
$5.50 \times 10^{-1}$ & $18558$ &
$9.97 \times 10^{-1}$ & $14048$ &
$5.79 \times 10^{-1}$ & $11278 $\\

& $100$ & $400$ &
$1.33 \times 10^{-5}$ & $18928$ &
$4.84 \times 10^{-1}$ & $169578$ &
$8.59 \times 10^{-3}$ & $18448$ &
$5.30 \times 10^{-1}$ & $67880 $\\

\hline
\multirow{2}{*}{$20000$} & $10$ & $50$ &
$1.95 \times 10^{-2}$ & $2589$ &
$7.30 \times 10^{-1}$ & $1475$ &
$ - $ & $ - $ & 
$7.20 \times 10^{-1}$ & $286 $\\

& $50$ & $200$ &
$1.49 \times 10^{-5}$ & $10364$ &
$6.30 \times 10^{-1}$ & $14588$ &
$ - $ & $ - $    &
$6.04 \times 10^{-1}$ & $29323$\\

\hline

$30000$ & $10$ & $50$ &
$1.62 \times 10^{-2}$ & $5767$ &
$7.74 \times 10^{-1}$ & $2437$ &
$ - $ & $ - $    &
$7.43 \times 10^{-1}$ & $308 $\\

\hline

\end{tabular}
\caption{ \footnotesize Numerical results for \optspace, SVT, FPCA and {\sc ADMiRA} for {\em hard} problems. $\eps = |E|/n$ is the number of observed entries per row/column.}\label{tab:hard}
\end{center}

\end{table*}
%\end{figure*}
\normalsize

\begin{table*}
\begin{center}
\footnotesize
\begin{tabular}{c c c | c c | c c | c c | c c }
\hline
\multirow{2}{*}{$n$} & \multirow{2}{*}{$r$} & \multirow{2}{*}{$\epsilon$} & \multicolumn{2}{|c|}{\optspace} & \multicolumn{2}{|c|}{SVT}& \multicolumn{2}{|c|}{FPCA}& \multicolumn{2}{|c}{{\sc ADMiRA}}\\
& & & rel. error & time(s) & rel. err. & time(s) & rel. err. & time(s) & rel. err. & time(s) \\
\hline
\multirow{3}{*}{$1000$} & $10$ & $120$ &
$1.18 \times 10^{-5}$ & $28$ &
$1.68 \times 10^{-5}$ & $40$ &
$5.20 \times 10^{-5}$ & $18$ &
$9.09 \times 10^{-4}$ & $52 $\\

& $50$ & $390$ &
$9.26 \times 10^{-6}$ & $212$ &
$1.62 \times 10^{-5}$ & $247$ &
$3.53 \times 10^{-6}$ & $106$  &
$3.62 \times 10^{-5}$ & $701 $\\

& $100$ & $570$ &
$1.49 \times 10^{-5}$ & $723$ &
$1.71 \times 10^{-5}$ & $694$ &
$1.92 \times 10^{-6}$ & $160$ &
$1.88 \times 10^{-5}$ & $2319 $\\
\hline

\multirow{3}{*}{$5000$} & $10$ & $120$ &
$1.51 \times 10^{-5}$ & $252$ &
$1.76 \times 10^{-5}$ & $112$ &
$1.69 \times 10^{-4}$ & $1083$ &
$4.68 \times 10^{-2}$ & $198 $\\

& $50$ & $500$ &
$1.16 \times 10^{-5}$ & $850$ &
$1.62 \times 10^{-5}$ & $1312$ &
$5.99 \times 10^{-5}$ & $1005$ &
$7.42 \times 10^{-3}$ & $92751 $\\

& $100$ & $800$ &
$8.39 \times 10^{-6}$ & $3714$ &
$1.73 \times 10^{-5}$ & $5432$ &
$3.32 \times 10^{-5}$ & $1953$ &
$4.42 \times 10^{-2}$ & $634028 $\\

\hline
\multirow{3}{*}{$10000$} & $10$ & $120$ &
$7.64 \times 10^{-6}$ & $632$ &
$1.75 \times 10^{-5}$ & $221$ &
$9.95 \times 10^{-1}$ & $13288$ &
$1.22 \times 10^{-1}$ & $442 $\\

& $50$ & $500$ &
$1.19 \times 10^{-5}$ & $2585$ &
$1.63 \times 10^{-5}$ & $2872$ &
$9.51 \times 10^{-5}$ & $7337$ &
$2.58 \times 10^{-2}$ & $186591 $\\

& $100$ & $800$ &
$1.46 \times 10^{-5}$ & $8514$ &
$1.76 \times 10^{-5}$ & $10962$ &
$6.90 \times 10^{-5}$ & $9426$ &
$9.66 \times 10^{-2}$ & $755082 $\\

\hline
\multirow{2}{*}{$20000$} & $10$ & $120$ &
$1.59 \times 10^{-5}$ & $1121$ &
$1.76 \times 10^{-5}$ & $461$ &
$ - $ & $ - $ &
$3.04 \times 10^{-1}$ & $181 $\\

& $50$ & $500$ &
$9.77 \times 10^{-6}$ & $4473$ &
$1.64 \times 10^{-5}$ & $6014$ &
$ - $ & $ - $    & 
$4.33 \times 10^{-2}$ & $346651$\\

\hline

$30000$ & $10$ & $120$ &
$1.56 \times 10^{-5}$ & $1925$ &
$1.80 \times 10^{-5}$ & $838$ &
$ - $ & $ - $    &
$4.19 \times 10^{-1}$ & $71 $\\

\hline
\end{tabular}
\end{center}
\caption{ \footnotesize Numerical results for \optspace, SVT, FPCA and {\sc ADMiRA} for {\em easy} problems. $\eps = |E|/n$ is the number of observed entries per row/column.}\label{tab:easy}

\end{table*}
\normalsize

%
%======================================================
%
%\subsubsection{Ill-conditioned matrices}
%\label{sec:result_incremental}

To add robustness in the case when the condition number of the matrix $M$ is high, 
we introduced a novel modification to {\sc OptSpace} 
in Section \ref{sec:IncrementalOptSpace}.
To illustrate the robustness of this {\sc Incremental OptSpace}, 
in Table \ref{tab:highcond}, we compare the results of exact matrix completion 
%obtained by the {\sc Incremental OptSpace} algorithm to that of other algorithms
for different values of $\kappa$.
Here, $\kappa$ denotes the condition number of 
the randomly generated matrix $M$ used in the simulation.
For this simulation with ill-conditioned matrices, 
we use $n \times n$ random matrices generated as follows. 
For fixed $n=1000$, let $\widetilde{U}\in\R^{n\times r}$ and $\widetilde{V}\in\R^{n\times r}$ be 
the orthonormal basis for the space spanned 
by the columns of $U$ and $V$ respectively. 
Also, let $D$ be an $r \times r$ diagonal 
matrix with its diagonal entries linearly 
spaces between $n$ and $n/\kappa$. 
Then the matrix $M$ is formed as 
$M = \widetilde{U} D \widetilde{V}^T$. 
We use $\delta_{\rm tol} = 10^{-5}$ as the stopping criterion.
Table \ref{tab:highcond} shows that {\sc Incremental OptSpace} improves significantly over 
{\sc OptSpace} and achieves results comparable to the other algorithms.

\begin{table*}
\addtolength{\oddsidemargin}{-1cm}
\addtolength{\evensidemargin}{-1cm}
%\begin{center}

\scriptsize
\begin{tabular}{c c|c c|c c|c c|c c|c c}
\hline
\multirow{2}{*}{$\kappa$} & \multirow{2}{*}{$r$} & \multicolumn{2}{|c|}{{\sc OptSpace}}  & \multicolumn{2}{|c|}{{\sc Inc. OptSpace}} & \multicolumn{2}{|c|}{SVT}& \multicolumn{2}{|c|}{FPCA} & \multicolumn{2}{|c}{{\sc ADMiRA}}\\
& & rel. error & time & rel. error & time(s) & rel. err. & time(s) & rel. err. & time(s) & rel. err. & time(s)  \\
\hline
\multirow{3}{*}{$1$} & $10$ &
$8.56 \times 10^{-6}$ & $20$ &
$8.66 \times 10^{-6}$ & $19$ &
$1.70 \times 10^{-5}$ & $55$ &
$5.32 \times 10^{-5}$ & $22$ &
$1.57 \times 10^{-5}$ & $242$ \\

& $50$  &
$1.16 \times 10^{-5}$ & $78$ &
$1.09 \times 10^{-5}$ & $832$ &
$1.64 \times 10^{-5}$ & $628$ &
$2.97 \times 10^{-6}$ & $115$ &
$1.60 \times 10^{-5}$ & $1252 $\\

& $100$ &
$7.05 \times 10^{-6}$ & $401$ &
$7.37 \times 10^{-6}$ & $4605$ &
$1.78 \times 10^{-5}$ & $2574$ &
$1.88 \times 10^{-6}$ & $174$ &
$1.67 \times 10^{-5}$ & $3454 $\\
\hline

\multirow{3}{*}{$5$} & $10$  &
$1.08 \times 10^{-1}$ & $124$ &
$1.53 \times 10^{-5}$ & $70$ &
$1.53 \times 10^{-5}$ & $72$ &
$5.53 \times 10^{-5}$ & $21$ &
$1.56 \times 10^{-5}$ & $234 $\\

& $50$  &
$1.10 \times 10^{-1}$ & $1591$ &
$1.30 \times 10^{-5}$ & $921$ &
$1.46 \times 10^{-5}$ & $639$ &
$1.08 \times 10^{-5}$ & $145$ &
$1.61 \times 10^{-5}$ & $1221 $\\

& $100$ &
$1.24 \times 10^{-1}$ & $5004$ &
$1.41 \times 10^{-5}$ & $5863$ &
$1.54 \times 10^{-5}$ & $1541$ &
$4.38 \times 10^{-6}$ & $664$ &
$1.68 \times 10^{-5}$ & $3450 $\\

\hline
\multirow{3}{*}{$10$} & $10$ & 
$1.09 \times 10^{-1}$ & $112$ &
$2.00 \times 10^{-1}$ & $238$ &
$1.47 \times 10^{-5}$ & $127$ &
$5.22 \times 10^{-5}$ & $21$ &
$1.55 \times 10^{-5}$ & $243 $\\

& $50$ & 
$1.04 \times 10^{-1}$ & $1410$ &
$1.32 \times 10^{-5}$ & $1593$ &
$1.36 \times 10^{-5}$ & $1018$ &
$1.42 \times 10^{-5}$ & $270$ &
$1.61 \times 10^{-5}$ & $1206 $\\

& $100$ &
$1.10 \times 10^{-1}$ & $4569$ &
$1.36 \times 10^{-5}$ & $9550$ &
$1.41 \times 10^{-5}$ & $2473$ &
$4.54 \times 10^{-6}$ & $996$ &
$1.65 \times 10^{-5}$ & $3426 $\\

\hline
\end{tabular}
%\end{center}
\caption{ \scriptsize Numerical results for {\sc OptSpace}, {\sc Incremental OptSpace}, SVT, FPCA and {\sc ADMiRA} for different condition numbers $\kappa$. $\eps = |E|/n$ depends only on $r$ and is the same as used in Table \ref{tab:noisy}. }\label{tab:highcond}

\addtolength{\oddsidemargin}{1cm}
\addtolength{\evensidemargin}{1cm}
\end{table*}

\normalsize
%
%======================================================
%
\subsection{Approximate matrix completion}
\label{sec:noisycompletion}

In this section we compare the performance of different algorithms 
for matrix completion with noisy observations.
As a metric, we use the relative root mean squared error defined as
\begin{eqnarray}
	{\rm RMSE} = \frac{1}{\sqrt{mn}}||M-\hM||_F \;. \label{eq:RMSE}
\end{eqnarray}

%
%======================================================
%
\subsubsection{Standard scenario}

For direct comparison we start with an example taken from \cite{CandesPlan}. 
In this example, $M$ is a square matrix of dimensions $n\times n$ and rank $r$ generated as 
$M = UV^T$ with fixed $n=600$. 
$U$ and $V$ are $n\times r$ matrices with 
each entry being sampled independently from a standard 
Gaussian distribution ${\cal N}(0,\sigma_s^2 = 20/ \sqrt{n})$.
As before, each entry is revealed independently 
with probability $\eps/n$. 
Each entry is corrupted by added noise matrix $Z$, so that the observation
for the index $(i,j)$ is $M_{ij}+Z_{ij}$.
Further, $Z$ has each entries drawn from 
i.i.d. standard Gaussian distribution ${\cal N}(0,1)$.
In the following we refer to this noise model as the standard scenario. 
We also refer to \cite{CandesPlan} for the data for the convex relaxation approach and 
the information theoretic lower bound. 

\begin{figure}
\begin{center}
\includegraphics[width=10cm]{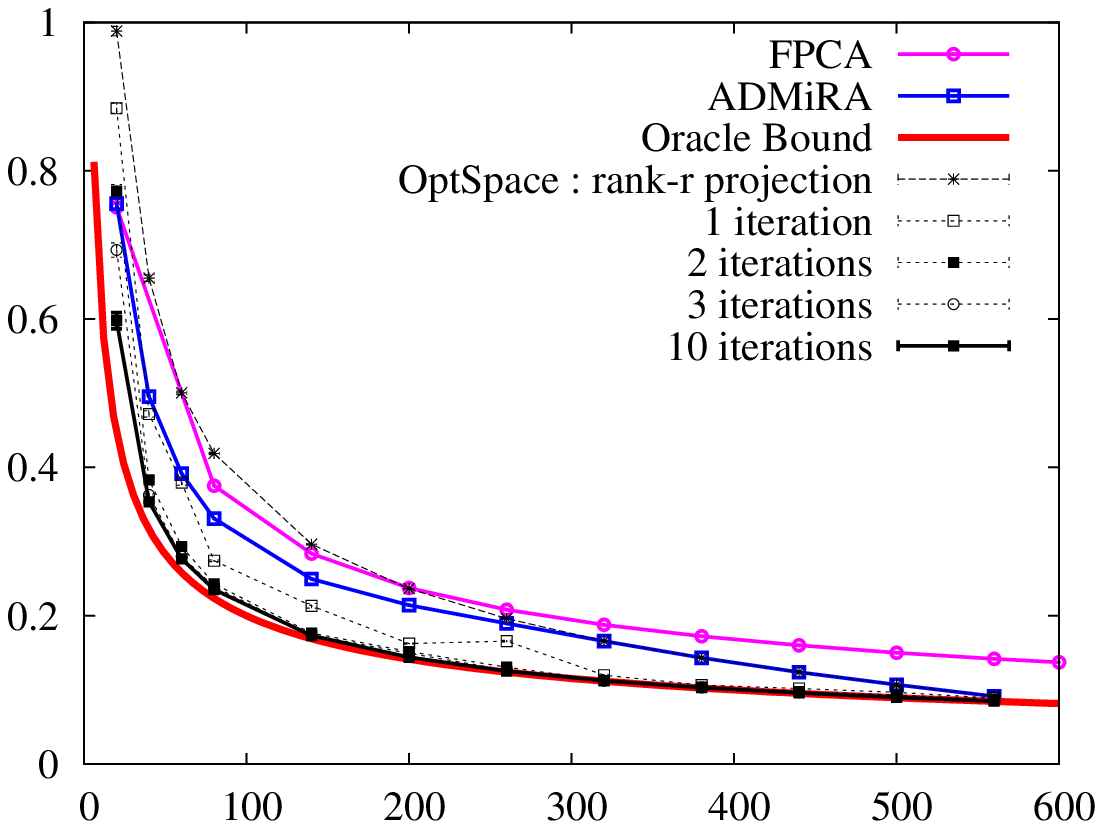}
%\hspace{-1.5cm}
%\input{noise_eps.tex}
\put(-134,-5){\scriptsize{$|E|/n$}}
\put(-265,88){\begin{sideways}\scriptsize{RMSE}\end{sideways}}
\caption{{\small Root mean squared error achieved by {\sc OptSpace} as a function of the observed entries $|E|$ and of the number of iterations in the {\sc Manifold Optimization} step. 
           $M$ is a rank-$2$ matrix with dimensions $m=n=600$.
		   The performances of the  convex relaxation approach and {\sc ADMiRA}, and an oracle lower bound are shown for comparison. }} \label{fig:noise_eps}
\end{center}
\end{figure}
\begin{figure}
\begin{center}
\includegraphics[width=10cm]{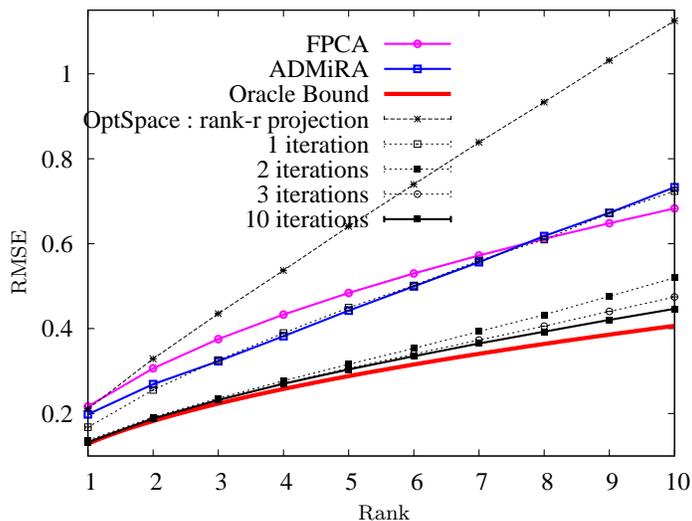}
%\hspace{-1.5cm}
%\input{noise_rank.tex}
\put(-134,-5){\scriptsize{Rank}}
\put(-265,80){\begin{sideways}\scriptsize{RMSE}\end{sideways}}
\caption{{\small The root mean square error of \optspace as a function of the rank $r$, and of the number of iterations in the {\sc Manifold Optimization} step. $M$ has dimensions $m=n=600$ and the number of observations is $|E|=72000$.  
The performances of the convex relaxation approach and {\sc ADMiRA}, and an oracle lower bound are shown for comparison. }} \label{fig:noise_rank}
\end{center}
\end{figure}

Figure \ref{fig:noise_eps} compares the average root mean squared error achieved by the different algorithms for a fixed rank $r=2$ as a function of $|E|/n$.  
After one iteration, for most values of $\eps$, \optspace has a smaller 
root mean square error than the convex relaxation approach and 
in about $10$ iterations, it becomes indistinguishable 
from the information theoretic lower bound. 
In Figure \ref{fig:noise_rank}, we compare the average root mean squared error obtained 
for a fixed sample size $\eps=120$ as a function of the rank. 
Again, for most values of $r$, after one iteration \optspace has a smaller 
root mean square error than the convex relaxation based algorithm. 
%For small values of $r$, \optspace becomes indistinguishable 
%from the information theoretic lower bound in about $10$ iterations. 
%For larger values of $r$, \optspace requires more samples to achieve the lower bound.
%The gap between \optspace and the lower bound at higher rank $r$ is expected, 
%since the number of samples required to get the same performance grows like $r$.

Table \ref{tab:noisy} illustrate how the performance 
changes with different noise power for fixed $n=1000$. 
We present the results of our experiments with different 
ranks and \emph{noise ratios} defined as 
\begin{eqnarray}
	N \equiv ||\cP_E(Z)||_F/||\cP_E(M)||_F\;, \label{eq:noiseratio}
\end{eqnarray}
as in \cite{CCS08}.

\begin{table*}
\begin{center}
\footnotesize
\begin{tabular}{c c c | c c | c c | c c | c c }
\hline
\multirow{2}{*}{$N$} & \multirow{2}{*}{$r$} & \multirow{2}{*}{$\epsilon$} & \multicolumn{2}{|c|}{\optspace} & \multicolumn{2}{|c|}{SVT}& \multicolumn{2}{|c|}{FPCA} & \multicolumn{2}{|c}{{\sc ADMiRA}}\\
& & & rel. error & time(s) & rel. err. & time(s) & rel. err. & time(s) & rel. err. & time(s)  \\
\hline
\multirow{3}{*}{$10^{-2}$} & $10$ & $120$ &
$4.47 \times 10^{-3}$ & $24$ &
$7.8 \times 10^{-3}$ & $11$ &
$5.48 \times 10^{-3}$ & $99$ &
$2.01 \times 10^{-2}$ & $35$ \\

& $50$ & $390$ &
$5.49 \times 10^{-3}$ & $149$ &
$9.5 \times 10^{-3}$ & $88$ &
$7.18 \times 10^{-3}$ & $805$ &
$1.83 \times 10^{-2}$ & $391 $\\

& $100$ & $570$ &
$6.39 \times 10^{-3}$ & $489$ &
$1.13 \times 10^{-2}$ & $216$ &
$1.08 \times 10^{-2}$ & $1111$ &
$1.63 \times 10^{-2}$ & $1424 $\\
\hline

\multirow{3}{*}{$10^{-1}$} & $10$ & $120$ &
$4.50 \times 10^{-2}$ & $23$ &
$0.72 \times 10^{-1}$ & $4$ &
$6.04 \times 10^{-2}$ & $140$ &
$1.18 \times 10^{-1}$ & $12 $\\

& $50$ & $390$ &
$5.52 \times 10^{-2}$ & $147$ &
$0.89 \times 10^{-1}$ & $33$ &
$7.77 \times 10^{-1}$ & $827$ &
$1.20 \times 10^{-1}$ & $139 $\\

& $100$ & $570$ &
$6.38 \times 10^{-2}$ & $484$ &
$1.01 \times 10^{-1}$ & $85$ &
$1.13 \times 10^{-1}$ & $1140$ &
$1.15 \times 10^{-1}$ & $572 $\\

\hline
\multirow{3}{*}{$1$} & $10$ & $120$ &
$4.86 \times 10^{-1}$ & $31$ &
$0.52 $ & $1$ &
$5.96 \times 10^{-1}$ & $141$ &
$5.38 \times 10^{-1}$ & $3 $\\

& $50$ & $390$ &
$6.33 \times 10^{-1}$ & $153$ &
$0.63 $ & $8$ &
$9.54 \times 10^{-1}$ & $1088$ &
$5.92 \times 10^{-1}$ & $47 $\\

& $100$ & $570$ &
$1.68 $ & $107$ &
$0.69 $ & $35$ &
$1.19 $ & $1582$ &
$6.69 \times 10^{-1}$ & $181 $\\

\hline
\end{tabular}
\end{center}
\caption{ \footnotesize Numerical results for \optspace, SVT, FPCA and {\sc ADMiRA} when the observations are corrupted by additive Gaussian noise with noise ratio $N$. $\eps = |E|/n$ is the number of observed entries per row/column. }\label{tab:noisy}

\end{table*}
\normalsize

Next, in the following series of examples, we illustrate 
how the performances change under different noise models. 
%We use $n\times n$ test matrices $M$ of rank $r$ generated as 
In the following, $M$ is a square matrix generated as $UV^T$ like above, 
but $U$ and $V$ now have each entry sampled independently from a standard 
Gaussian distribution ${\cal N}(0,1)$, unless specified otherwise.
As before, each entry is revealed independently with probability $\eps/n$
and the observation is corrupted by added noise matrix $Z$.
We compare the resulting RMSE of {FPCA}, {\sc ADMiRA} and {\sc OptSpace}
on this randomly generated matrices with noisy observations and missing entries.
Since {\sc ADMiRA} requires a target rank, we use the rank estimated using 
{\sc Rank Estimation} described in Section \ref{sec:rankestimation}.
For {FPCA} we choose $\mu=\sqrt{2np}\sigma$, where $p=|E|/n^2$ and $\sigma^2$
is the variance of each entry in $Z$. A convincing argument for 
this choice of $\mu$ is given in \cite{CandesPlan}. 

In the standard scenario, 
we typically make the following three assumptions on the noise matrix $Z$.
(1) The noise $Z_{ij}$ does not depend on the value of the matrix $M_{ij}$.
(2) The entries of $Z$, $\{Z_{ij}\}$, are independent.
(3) The distribution of each entries of $Z$ is Gaussian.
The matrix completion algorithms described in Section \ref{sec:model} 
are expected to be especially effective 
under this standard scenario for the following two reasons. 
First, the squared error objective function that the algorithms minimize is 
well suited for the Gaussian noise. Second, the independence of $Z_{ij}$'s 
ensure that the noise matrix is almost full rank and 
the singular values are evenly distributed. This implies that 
for a given noise power $||Z||_F$, 
the spectral norm $||Z||_2$ is much smaller than $||Z||_F$.
In the following, we fix $m=n=500$ and $r=4$, 
%but vary one or more assumptions on the noise 
%from the standard setting at a time to 
and study how the performance changes with different noise.
Each of the simulation results is averaged over 10 instances
and is shown with respect to two basic parameters,
the average number of revealed entries per row $|E|/n$ and 
the {noise ratios} $N$, defined as Eq.~(\ref{eq:noiseratio}).
%the signal-to-noise ratio, ${\rm SNR}= \E[||M||^2_F]/\E[||Z||^2_F] $.
%the square root of the expected noise power $NP=\sqrt{\E[||Z||^2_F/n^2]}$.

%
%======================================================
%
%\subsubsection{Standard scenario}
%{\bf Standard scenario.}
%
\begin{figure}
\begin{center}
\includegraphics[width=10cm]{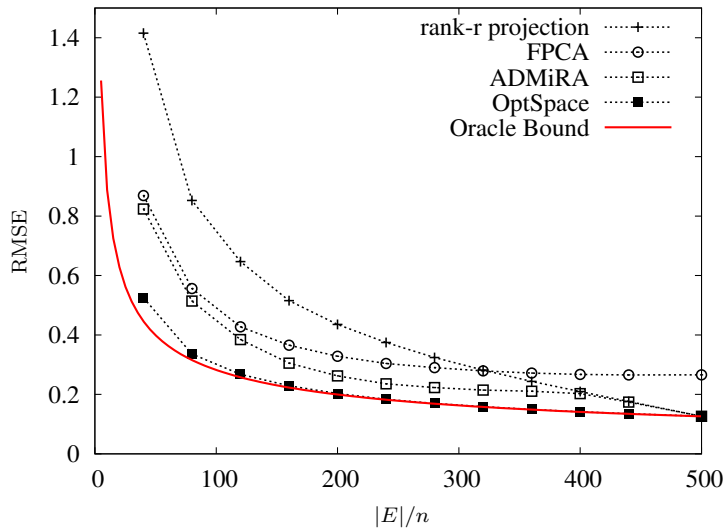}
\put(-140,-7){\scriptsize{$|E|/n$}}
\put(-277,87){\begin{sideways}\scriptsize{RMSE}\end{sideways}}\\
%\put(-114,0){\scriptsize{$|E|/n$}}
%\put(-220,80){\begin{sideways}\scriptsize{RMSE}\end{sideways}}\\
\hspace{-0.3cm}
%\includegraphics[width=8.3cm]{FIG/noise0_eps_t.eps}
%\put(-106,0){\scriptsize{$\eps$}}
%\put(-220,74){\begin{sideways}\scriptsize{seconds}\end{sideways}}\\
\caption{{\small The root mean squared error as a function of the average number of observed entries per row $|E|/n$ for fixed noise ratio $N=1/2$ under the standard scenario. }} \label{fig:noise0_eps}
\end{center}
\end{figure}

In this standard scenario, the noise $Z_{ij}$'s are distributed as 
i.i.d. Gaussian {\cal N}(0,$\sigma^2$). 
Note that the noise ratio is equal to $N=\sigma/2$.
%There is a basic trade-off between two metrics of interest:
The accuracy of the estimation is measured using RMSE.
We compare the resulting RMSE of {FPCA}, {\sc ADMiRA} and {\sc OptSpace}
to the RMSE of the oracle estimate, which is $\sigma\sqrt{(2nr-r^2)/|E|}$  \cite{CandesPlan}.

%and the computation complexity is measured by the running time in seconds.

%In order to interpret the simulation results, 
%they are compared to the RMSE achieved by the oracle
%and a simple rank-$r$ projection algorithm defined as Eq.~(\ref{eq:projection}).
%The rank-$r$ projection algorithm simply computes $\cP_r(N^E)$.
%the SVD of $N^E$
%and outputs the rank-$r$ principal component normalized by $n/\eps$, defined as Eq.~(\ref{eq:projection}).
%The oracle has prior knowledge of  
%the linear subspace spanned by $\{UX^T+YV^T:X\in\R^{m\times r},Y\in\R^{n\times r}\}$, and 
%the RMSE of the oracle estimate is $\sigma\sqrt{(2nr-r^2)/n\eps}$  \cite{CandesPlan}.

Figure \ref{fig:noise0_eps} shows the performance for each of 
the algorithms with respect to $|E|/n$ under the standard scenario for fixed $N=1/2$. 
For most values of $|E|$, the simple rank-$r$ projection has the worst performance. 
However, when all the entries are revealed and the noise is i.i.d. Gaussian, 
the simple rank-$r$ projection coincides with the oracle bound, 
which in this simulation corresponds to the value $|E|/n=500$. 
Note that the behavior of the performance curves of FPCA, {\sc ADMiRA}, 
and {\sc OptSpace} with respect to $|E|$ is similar to the oracle bound, 
which is proportional to $1/\sqrt{|E|}$. 

Among the three algorithms, FPCA has the largest RMSE, 
and {\sc OptSpace} is very close to the oracle bound for all values of $|E|$. 
Note that when all the values are revealed, {\sc ADMiRA} is an efficient way 
of implementing rank-$r$ projection, and the performances are expected to be similar. 
This is confirmed by the observation that for $|E|/n \geq 400$ the two curves are 
almost identical. One of the reasons why the RMSE of FPCA 
does not decrease with $|E|$ for large values of $|E|$ is that 
FPCA overestimates the rank and returns estimated matrices 
with rank much higher than $r$, whereas the rank estimation algorithm
 used for {\sc ADMiRA} and {\sc OptSpace} always returned the 
 correct rank $r$ for $|E|/n \geq 80$. 

%The second figure in Figure \ref{fig:noise0_eps} shows 
%the average running time of the algorithms with respect to $\eps$.
%Note that due to the large difference between the running time of three algorithms, 
%the time is displayed in log scale.
%For most of the simulations, {\sc ADMiRA} had shortest running time and {FPCA} the 
%longest, and the gap was noticeably large as clearly shown in the figure.
%For FPCA and {\sc OptSpace}, the computation time increased with $\eps$, whereas 
%{\sc ADMiRA} had relatively stable computation time independent of $\eps$.

%
\begin{figure}
\begin{center}
\includegraphics[width=10cm]{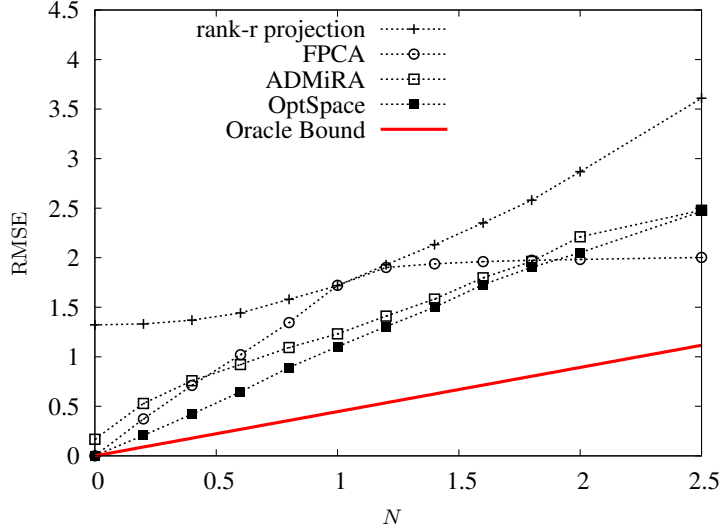}
\put(-137,-7){\scriptsize{$N$}}
\put(-277,87){\begin{sideways}\scriptsize{RMSE}\end{sideways}}\\
%\put(-106,0){\scriptsize{$N$}}
%\put(-220,80){\begin{sideways}\scriptsize{RMSE}\end{sideways}}\\
\hspace{-0.1cm}
%\includegraphics[width=8.1cm]{FIG/noise0_eps40sig_t.eps}
%\put(-121,0){\scriptsize{$1/\sqrt{\rm SNR}$}}
%\put(-220,74){\begin{sideways}\scriptsize{seconds}\end{sideways}}\\
\caption{{\small The root mean squared error as a function of noise ratio $N$ for fixed $|E|/n=40$ under the standard scenario.}} \label{fig:noise0_sig}
\end{center}
\end{figure}

Figure \ref{fig:noise0_sig} show the performance for each of 
the algorithms against the noise ratio $N$
within the standard scenario for fixed $|E|/n=40$. 
The behavior of the performance curves of {\sc ADMiRA} and {\sc OptSpace} are
similar to the oracle bound which is linear in $\sigma$ 
which, in the standard scenario, is equal to $2N$. 
The performance of the rank-$r$ projection algorithm is 
determined by two factors. One is the added noise which is linear in $N$
and the other is the error caused by the erased entries which is 
constant independent of $N$. These two factors add up, whence the 
performance curve of the rank-$r$ projection follows. 
The reason the RMSE of FPCA does not decrease with SNR
for values of SNR less than $1$ is not that the estimates are good but
rather the estimated entries gets very small and the resulting RMSE is close to 
$\sqrt{\E[||M||_F^2/n^2]}$, which is $2$ in this simulation, regardless of the noise power.
When there is no noise, which corresponds to the value $N=0$, 
FPCA and {\sc OptSpace} both recover the original matrix correctly 
for this chosen value of $|E|/n=40$.
%In the large noise regime 
%norm( M,'fro')/n    1.9748
%FPCA    0.4302
%OptSpace    2.0625
%ADMiRA    1.9930
%SVD    3.5908
%Hence FPCA seems to be pretty independent of the noise, but it does not imply that ...
%It only controls overall deviation from zero....
%For all three algorithms, the computation time is larger for smaller noise,
%and the reason is that it takes more iterations until the stopping criterion is met.
%Also, for most of the simulations with different SNR, 
%{\sc ADMiRA} had shortest running time and 
%{FPCA} the longest.

%
%======================================================
%
\subsubsection{Multiplicative Gaussian noise}
%{\bf Multiplicative Gaussian noise.}
%
\begin{figure}
\begin{center}
\includegraphics[width=10cm]{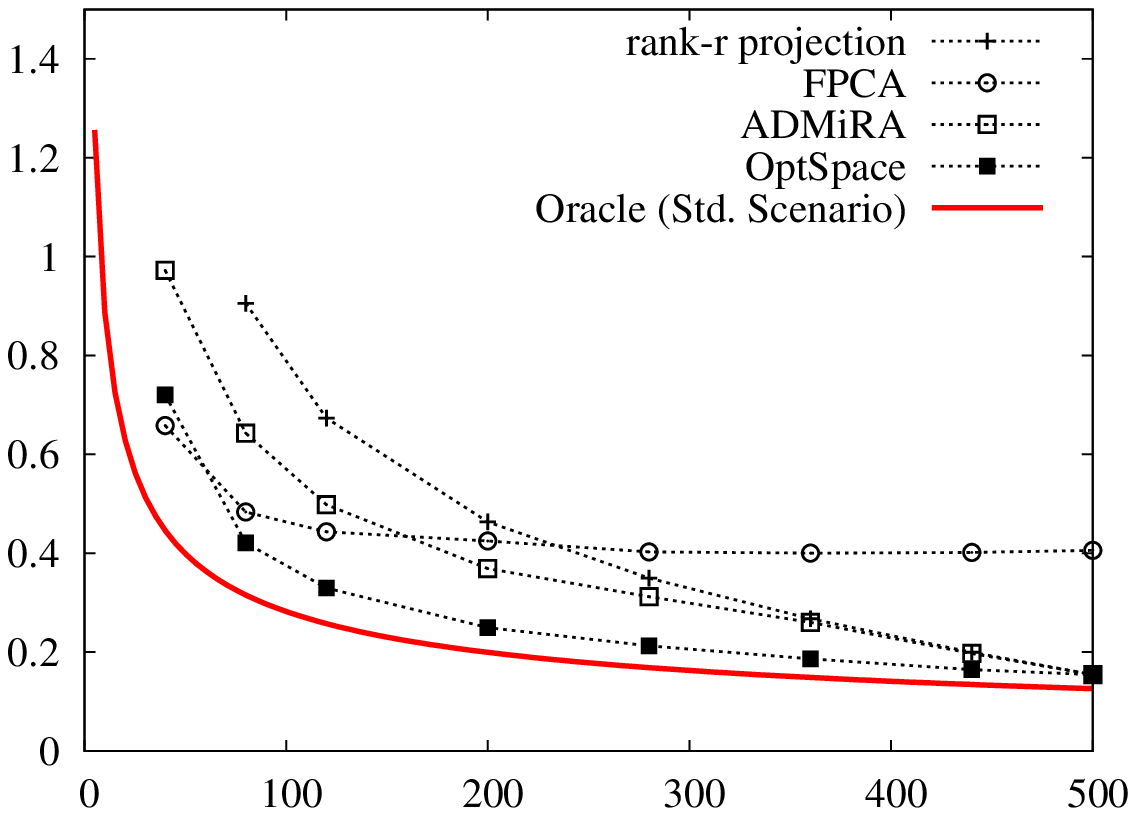}
\put(-140,-7){\scriptsize{$|E|/n$}}
\put(-277,87){\begin{sideways}\scriptsize{RMSE}\end{sideways}}
%\put(-114,0){\scriptsize{$|E|/n$}}
%\put(-220,80){\begin{sideways}\scriptsize{RMSE}\end{sideways}}
\hspace{-0.3cm}
\caption{{\small The root mean squared error as a function of the average number of observed entries per row $|E|/n$ for fixed noise ratio $N=1/2$ within the multiplicative noise model. }} \label{fig:noise_multi}
\end{center}
\end{figure}

In sensor network localization \cite{WZBY06}, where the entries of the matrix corresponds to 
the pair-wise distances between the sensors, the observation noise is oftentimes 
assumed to be multiplicative. In formulae, $Z_{ij}=\xi_{ij}M_{ij}$, 
where $\xi_{ij}$'s are distributed as i.i.d. Gaussian with zero mean.
The variance of $\xi_{ij}$'s are chosen to be $1/r$ 
so that the resulting noise ratio is $N=1/2$.
Note that in this case, $Z_{ij}$'s are mutually dependent through $M_{ij}$'s and 
the values of the noise also depend on the value of the matrix entry $M_{ij}$.

Figure \ref{fig:noise_multi} shows the RMSE with respect to 
$|E|/n$ under multiplicative Gaussian noise. The RMSE of the rank-$r$ projection for $|E|/n=40$
is larger than $1.5$ and is omitted in the figure.
The bottommost line corresponds to the oracle performance under standard scenario,
and is displayed here, and all of the following figures, to serve as a reference for comparison.
The main difference with respect to Figure \ref{fig:noise0_eps} is that 
most of the performance curves are larger under multiplicative noise.
For the same value of the noise ratio $N$, 
it is more difficult to distinguish the noise from the original matrix,
since the noise is now correlated with the matrix $M$.

%
%======================================================
%
\subsubsection{Outliers}
%{\bf Outliers.}

%
\begin{figure}
\begin{center}
\includegraphics[width=8cm]{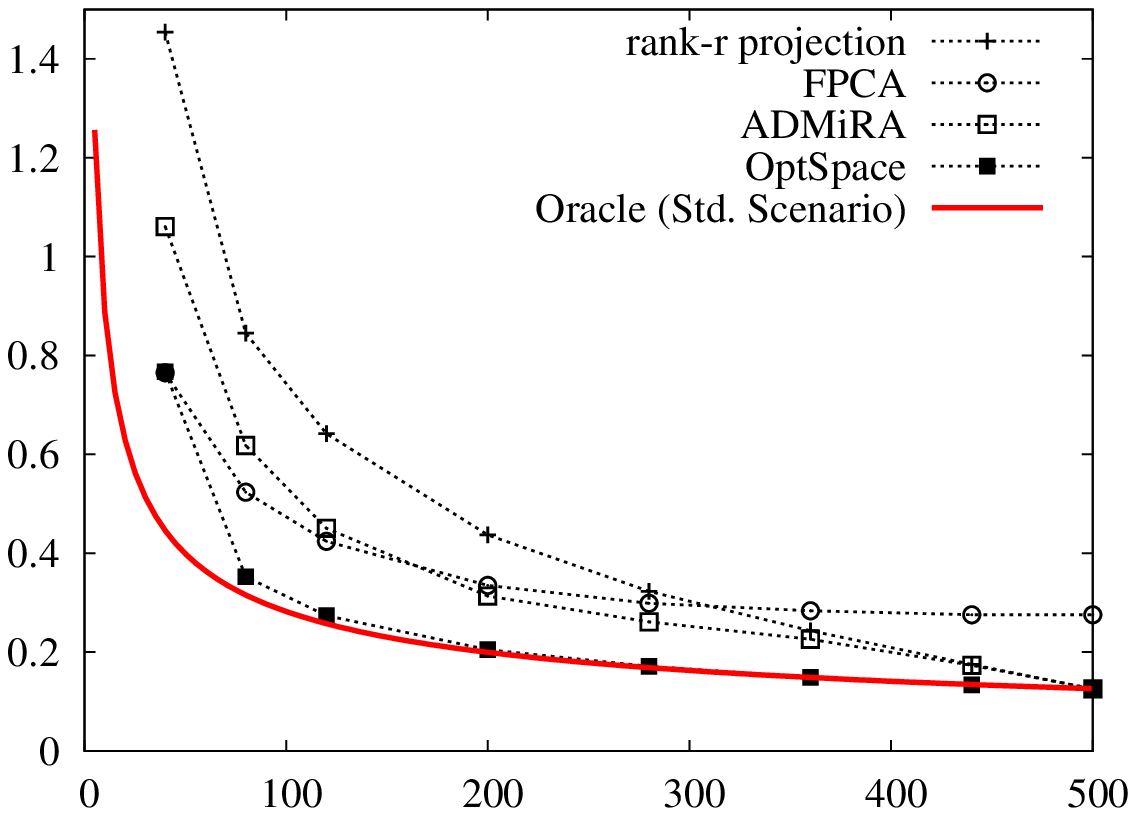}
%\put(-140,-7){\scriptsize{$|E|/n$}}
%\put(-277,87){\begin{sideways}\scriptsize{RMSE}\end{sideways}}
\put(-114,0){\scriptsize{$|E|/n$}}
\put(-220,80){\begin{sideways}\scriptsize{RMSE}\end{sideways}}
\includegraphics[width=8cm]{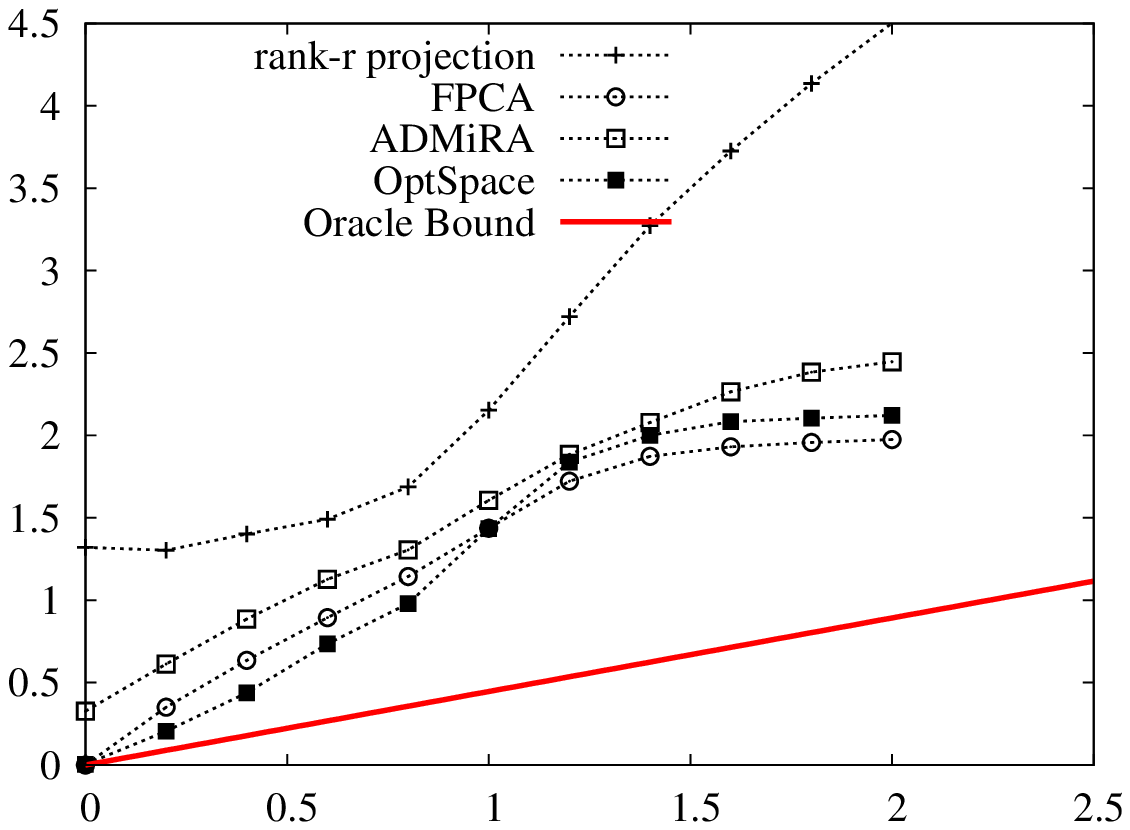}
%\put(-137,-7){\scriptsize{$N$}}
%\put(-277,87){\begin{sideways}\scriptsize{RMSE}\end{sideways}}
\put(-106,0){\scriptsize{$N$}}
\put(-220,80){\begin{sideways}\scriptsize{RMSE}\end{sideways}}
%\hspace{-0.3cm}
\caption{{\small The root mean squared error as a function of the average number of observed entries per row $|E|/n$ for fixed noise ratio $N=1/2$ with outliers (left) and the RMSE as a function of the noise ratio $N$ for fixed $|E|/n=40$ with outliers (right).}} \label{fig:noise_outlier}
\end{center}
\end{figure}
\vspace{-0.2cm}
In structure from motion \cite{SFM}, the entries of the matrix corresponds to 
the position of points of interest in $2$-dimensional images captured by 
cameras in different angles and locations. 
However, due to failures in the feature extraction algorithm, 
some of the observed positions are corrupted by large noise 
where as most of the observations are noise free.
To account for such outliers, we use the following model. 
\begin{eqnarray*}
  Z_{ij} = \left\{
              \begin{array}{rl}
              a  & \text{with probability $1/200$ } \, ,\\
              -a & \text{w.p. $1/200$ } \, ,\\
              0       & \text{w.p. $99/100$}\, .
              \end{array} \right. 
\end{eqnarray*}
The value of $a$ is chosen according to the target noise ratio $N=a/20$.
The noise is independent of the matrix entries and 
$Z_{ij}$'s are mutually independent, but the distribution is now non-Gaussian. 

Figure \ref{fig:noise_outlier} shows the performance of 
the algorithms with respect to $|E|/n$ and the noise ratio $N$ with outliers.
Comparing the first figure to Figure \ref{fig:noise0_eps}, we can see that
the performance for large value of $|E|$ is less affected by outliers compared to 
the small values of $|E|$. The second figure clearly shows how the performance degrades for non-Gaussian noise when the number of samples is small.  
The algorithms minimize the squared error $||\cP_E(X)-\cP_E(N)||_F^2$ as in (\ref{P3}) and (\ref{P4}).
For outliers, a suitable algorithm would be to 
minimize the $l_1$-norm of the errors instead of the $l_2$-norm \cite{CSPW09,WGRM09}. 
Hence, for this simulation with outliers, we can see that the performance of the rank-$r$ 
projection, {\sc ADMiRA} and {\sc OptSpace} is worse than the Gaussian noise case.
However, the performance of FPCA is almost the same as in the standard scenario. 
%slightly increased in simulation with outliers, 
%and this is related to the fact that FPCA tends to give estimations which have smaller values.

%
%======================================================
%
\subsubsection{Quantization noise}
%{\bf Quantization noise.}

%
\begin{figure}
\begin{center}
\includegraphics[width=10cm]{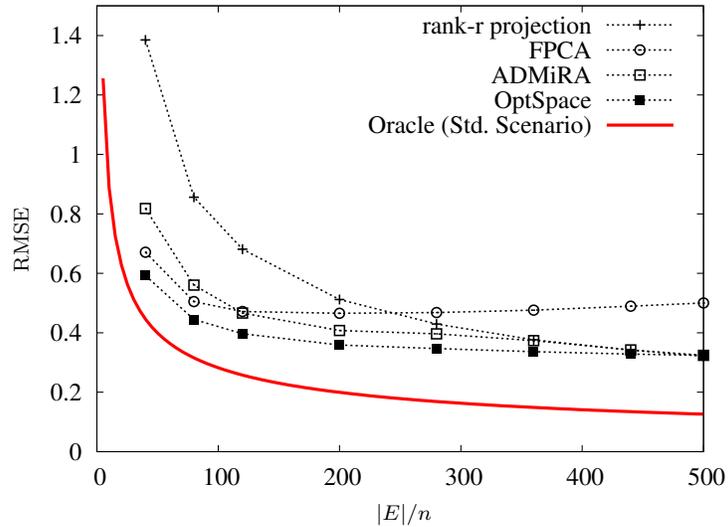}
\put(-140,-7){\scriptsize{$|E|/n$}}
\put(-277,87){\begin{sideways}\scriptsize{RMSE}\end{sideways}}\\
%\put(-114,0){\scriptsize{$|E|/n$}}
%\put(-220,80){\begin{sideways}\scriptsize{RMSE}\end{sideways}}\\
\caption{{\small The root mean squared error as a function of the average number of observed entries per row $|E|/n$ for fixed noise ratio $N=1/2$ with quantization.}} \label{fig:noise_quantization}
\end{center}
\end{figure}

One common model for noise is the quantization noise. 
For a regular quantization, we choose a parameter $a$ and 
quantize the matrix entries to the nearest value in  $\{\ldots$, $-a/2$, $a/2$, $3a/2$, $5a/2$, $\ldots\}$.
The parameter $a$ is chosen carefully such that the resulting noise ratio is $1/2$.
The performance for this quantization 
is expected to be worse than the multiplicative noise case.
The reason is that now the noise is deterministic and completely depends on the matrix entries $M_{ij}$, 
whereas in the multiplicative noise model it was random.

Figure \ref{fig:noise_quantization} shows the performance against $|E|/n$
within quantization noise. The overall behavior of the performance curves is similar to Figure
\ref{fig:noise0_eps}, but most of the curves are shifted up. Note that the bottommost line 
is the oracle performance in the standard scenario which is the same in all the figures.
Compared to Figure \ref{fig:noise_multi}, for the same value of $N=1/2$,  
quantization is much more detrimental than the multiplicative noise as expected.

%
%======================================================
%
\subsubsection{Ill conditioned matrices}
%{\bf Ill conditioned matrices.}

%
\begin{figure}
\begin{center}
\includegraphics[width=10cm]{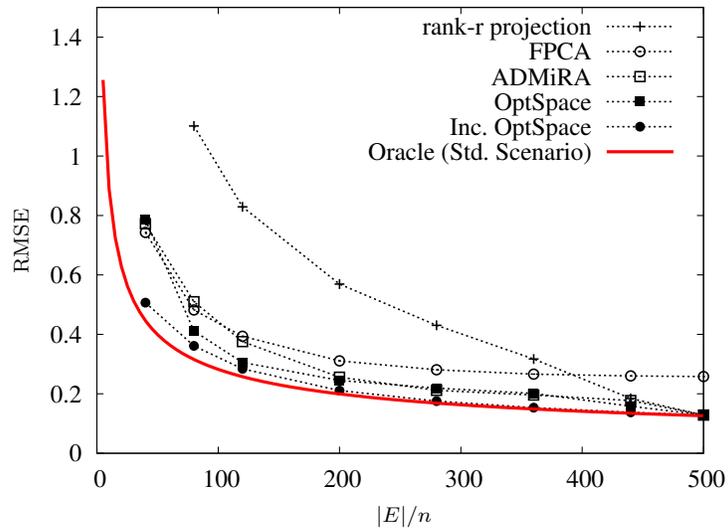}
\put(-140,-7){\scriptsize{$|E|/n$}}
\put(-277,87){\begin{sideways}\scriptsize{RMSE}\end{sideways}}\\
%\put(-114,0){\scriptsize{$|E|/n$}}
%\put(-220,80){\begin{sideways}\scriptsize{RMSE}\end{sideways}}\\
\caption{{\small The root mean squared error as a function of the average number of observed entries per row $|E|/n$ for fixed $N=1/2$ with ill-conditioned matrices.}} \label{fig:noise_ill}
\end{center}
\end{figure}

In this simulation, we look at how the performance degrades under the standard scenario 
if the matrix $M$ is ill-conditioned. $M$ is generated as 
$M=\sqrt{4/166}\, U \,\diag([1,4,7,10])\, V^T$, where $U$ and $V$ are generated 
as in the standard scenario. The resulting matrix has a condition number close to $10$
and the normalization constant $\sqrt{4/166}$ is chosen such that $\E[||M||_F]$ 
is the same as in the standard case. 
Figure \ref{fig:noise_ill} shows the performance as a function of $|E|/n$
with ill-conditioned matrix $M$. 
The performance of {\sc OptSpace} is similar to that of {\sc ADMiRA} for many values of $|E|$. 
However, similar to the noiseless simulation results, {\sc Incremental OptSpace} 
achieves a better performance in this case of ill-conditioned matrix.
%The {\sc Incremental OptSpace} algorithm starts from finding a rank-$1$ approximation
%from $N^E$ and incrementally finds higher rank approximations and 
%has more robust performance when $M$ is ill-conditioned, but is computationally more expensive.

%
%======================================================
%

\section{Numerical results with real data matrices}
\label{sec:real}

In this section, we consider low-rank matrix completion problems in the context of recommender systems, 
based on two real data sets: 
the Jester joke data set \cite{jesterjokes} and the Movielens data set \cite{movielens}. 
The Jester joke data set contains $4.1\times 10^6$ ratings for 
100 jokes from 73,421 users. \footnote{The dataset is available at {http://www.ieor.berkeley.edu/$\sim$goldberg/jester-data/}} Since the number of users 
is large compared to the number of jokes, we randomly select 
$n_u\in\{100,1000,2000,4000\}$ users for comparison purposes.
As in \cite{FPCA}, we randomly choose two ratings for each user as a test set, 
and this test set, which we denote by $T$, 
is used in computing the prediction error in Normalized Mean Absolute Error (NMAE). 
The Mean Absolute Error (MAE) is defined as in \cite{FPCA,eigentaste}. 
\begin{eqnarray*}
 MAE = \frac{1}{|T|}\sum_{(i,j)\in T} |M_{ij}-\hM_{ij}| \;,
\end{eqnarray*}
where $M_{ij}$ is the original rating in the data set and $\hM_{ij}$ 
is the predicted rating for user $i$ and item $j$.
The Normalized Mean Absolute Error (NMAE) is defined as 
\begin{eqnarray*}
 NMAE = \frac{MAE}{M_{\rm max}-M_{\rm min}} \;,
\end{eqnarray*}
where $M_{\rm max}$ and $M_{\rm min}$ are upper and lower bounds for the ratings. 
In the case of Jester joke, all the ratings are in $[-10,10]$ which implies that 
$M_{\rm max}=10$ and $M_{\rm min}=-10$. 

The numerical results for Jester joke data set using {\sc Incremental OptSpace}, 
{\sc FPCA} and {\sc ADMiRA} are presented in the first four columns of Table \ref{tab:jesterjoke}.
In the table, $rank$ indicates the rank used to estimate the matrix and $time$ 
is the running time of each matrix completion algorithm. To get an idea of how good the predictions are, 
consider the case where each missing entries is predicted with a random number 
drawn uniformly at random in $[-10,10]$ and the actual rating is also 
a random number with same distribution.
After a simple computation, we can see that 
the resulting NMAE of the random prediction is 0.333.
As another comparison, for the same data set with $n_u=18000$, 
simple nearest neighbor algorithm and {\sc Eigentaste} both 
yield NMAE of 0.187 \cite{eigentaste}. The NMAE of \incoptspace is lower than 
these simple algorithms even for $n_u=100$ and tends to decrease with $n_u$.

\begin{table*}
\begin{center}
\footnotesize
\begin{tabular}{c c c c | c c | c c | c c }
\hline
\multirow{2}{*}{$n_u$} & \multirow{2}{*}{m} &\multirow{2}{*}{samples} & \multirow{2}{*}{rank} & \multicolumn{2}{|c|}{\incoptspace} & \multicolumn{2}{|c|}{FPCA} & \multicolumn{2}{|c}{{\sc ADMiRA}}\\
& & & & NMAE & time(s) & NMAE & time(s) & NMAE & time(s)  \\
\hline
%\multirow{4}{*}{Jester joke}&
$100$ & $100$ & $7484$ & $2$ &
$0.17674 $ & $0.1$ &
$0.20386 $ & $25$ &
$0.18194 $ & $0.3$ \\

$1000$ & $100$ & $73626$ & $9$ &
$0.15835 $ & $11$ &
$0.16114 $ & $111$ &
$0.16194 $ & $0.5$ \\

$2000$ & $100$ & $146700$ & $9$ &
$0.15747 $ & $26$ &
$0.16101 $ & $243$ &
$0.16286 $ & $0.9 $\\

$4000$ & $100$ & $290473$ & $9$ &
$0.15918 $ & $56$ &
$0.16291 $ & $512$ &
$0.16317 $ & $2 $\\
\hline
%Movielens &
$943$ & $1682$ & $80000$ & $10$ &
$0.18638 $ & $213$ &
$0.19018 $ & $753$ &
$0.24276 $ & $5$ \\
\hline

\end{tabular}
\end{center} 
\caption{ \footnotesize Numerical results for the Jester joke data set, 
           where the number of jokes m is fixed at 100 (top four rows), 
           and for the Movielens data set with $943$ users and $1682$ movies (last row).
        }\label{tab:jesterjoke}
\end{table*}
\normalsize

Numerical simulation results on the Movielens data set
is also shown in the last row of the above table. 
The data set contains $100,000$ ratings for $1,682$ movies from $942$ users.\footnote{The dataset is available at {http://www.grouplens.org/node/73}} We use 
$80,000$ randomly chosen ratings to estimate the $20,000$ ratings in the test set, 
which is called $u1.base$ and $u1.test$, respectively, in the movielens data set. 
In the last column of Table \ref{tab:jesterjoke}, we compare the resulting NMAE using \incoptspace, {\sc FPCA} and {\sc ADMiRA}.

\begin{figure}
\begin{center}
%\hspace{-1.5cm}
\includegraphics[width=8cm]{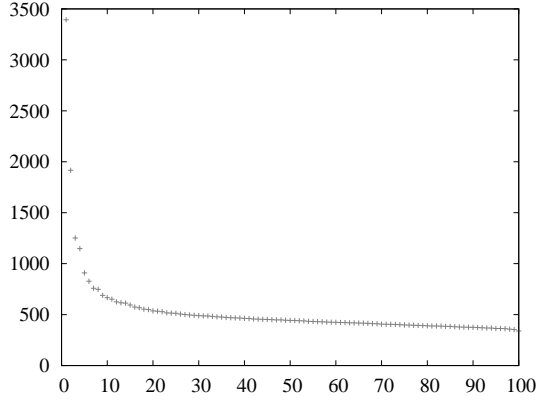}
\caption{{\small Distribution of the singular values of the complete sub matrix in the Jester joke data set. }} \label{fig:jesterjoke}
\end{center}
\end{figure}

Next, to get some insight on the structure of real data, 
we look at a complete sub matrix where all the entries are known. 
With Jester joke data set, we deleted all users containing 
missing entries, and generated a complete matrix $M$ 
with $14,116$ users and $100$ jokes. 
The distribution of 
the singular values of $M$ is shown in Figure \ref{fig:jesterjoke}. 
We must point out that this rating matrix is not low-rank or even approximately low-rank, 
although it is common to make such assumptions. This is one of  the 
difficulties in dealing with real data. The other aspect is that 
the samples are not drawn uniformly at random as commonly assumed in \cite{CandesTaoMatrix,KOM09}.

\begin{figure}
\begin{center}
%\hspace{-1.5cm}
\includegraphics[width=8cm]{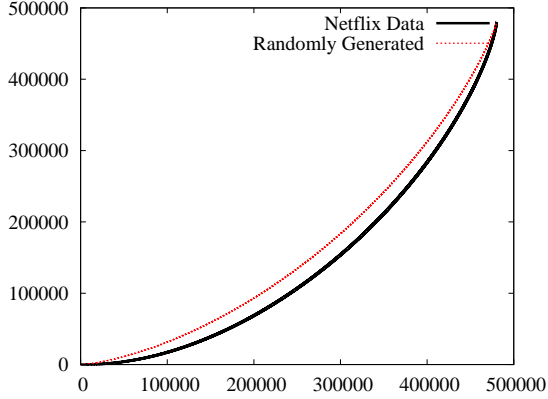}
\caption{{\small Cumulative sum of the sorted row norms of the factor corresponding to users.}} \label{fig:usersinc}
\end{center}
\end{figure}

\begin{figure}
\begin{center}
%\hspace{-1.5cm}
\includegraphics[width=8cm]{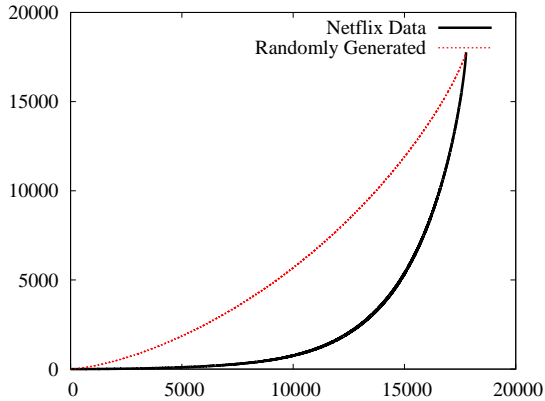}
\caption{{\small Cumulative sum of the sorted row norms of the factor corresponding to movies.}} \label{fig:movsinc}
\end{center}
\end{figure}

Finally we test the incoherence assumption for the Netflix 
dataset \cite{Net06} in Figure \ref{fig:usersinc} and Figure \ref{fig:movsinc}. 
For a $m\times n$ matrix whose singular value decomposition is given by $M=U\Sigma V^T$,
$M$ is said to be $(\mu_0,\mu_1)$-incoherent if it satisfies the following properties \cite{CaR08}:
\begin{itemize}
\item[{\bf A1.}] There exists a constant 
$\mu_0>0$ such that for all  $i\in [m]$, $j\in [n]$ we have 
           $\sum_{k=1}^{r}{U_{i,k}^2} \le \mu_0 r/m$,
           $\sum_{k=1}^{r}{V_{i,k}^2} \le \mu_0 r/n$.
\item[{\bf A2.}] There exists $\mu_1$ such that $|\sum_{k=1}^{r}{U_{i,k}V_{j,k}}|\leq\mu_1\sqrt{r/mn}$.
\end{itemize}
To check if {\bf A1} holds for the Netflix movie ratings matrix, 
we run OptSpace on the Netflix dataset and plot cumulative sum of 
the sorted row norms of the left and right factors defined as follows. 
Let the output of OptSpace be $X \in \reals^{m \times r}$, $Y \in \reals^{n\times r}$ and $S \in \reals^{r\times r}$.
Here $m = 480,189$ is the number of users and $n = 17,770$ is the number of movies.
For the target rank we used $r = 5$. 
Let $x_i = \frac{m}{r} ||X^{(i)}||^2$ and $y_i = \frac{n}{r} ||Y^{(i)}||^2$ where
$X^{(i)}$ and $Y^{(i)}$ denote the $i$th row of the left factor $X$ and the right factor $Y$ respectively.
Define a permutation $\pi_l:[m]\rightarrow [m]$ which sorts $x_i$'s in a non-decreasing order such that 
$x_{\pi_l(1)} \le x_{\pi_l(2)} \le \ldots \le x_{\pi_l(m)}$.
Here, we used the standard combinatorics notation $[k]=\{1,2,\ldots,k\}$ for an integer $k$.
Similarly, we can define $\pi_r:[n]\rightarrow [n]$ for $y_i$'s. 

In Figure \ref{fig:usersinc}, we plot $\sum_{i=1}^{k} x_i$  vs. $k$. 
For comparison, we also plot the corresponding
results for a randomly generated matrix $X_G$. 
Generate $U \in \reals^{m\times r}$ by sampling its entries $U_{ij}$
independently and distributed as ${\cal N}(0,1)$ and 
let $X_G$ be the left singular vectors of $U$. 
Since $x_i$'s are scaled by $m/r$, when $k=m$ we have 
$\sum_{i=1}^{m} x_i=m$. This is also true for the random matrix $X_G$.
Figure \ref{fig:movsinc} shows the corresponding plots for $Y$. 
For a given matrix, if {\bf A1} holds with a small $\mu_0$ then 
the corresponding curve would be close to a straight line. 
The curvature in the plots is indicative of the disparity among
the row weigths of the factors. 
We can see that a randomly generated matrix would satisfy 
{\bf A1} with a smaller value of $\mu_0$ compared to the movie ratings matrix,
hence can be said to be more incoherent. 
The factor corresponding to movies has a larger disparity
than the factor corresponding to users,
and hence challenges the incoherence assumption.

%
%======================================================
%

%\section{Extensions to real data}

%1. What is the difficulties with dealing with REAL DATA.\\
%2. Some Smart Tricks we used. \\
%-Subtracting Mean \\
%-Rescaling \\
%3. How should we interpret real data? : low rank + noise? \\
%4. Residual Noise Statistics.\\

%\newpage
\appendix

%
% =============================================================================
%

\section{Proof of Proposition \ref{pro:rank}} \label{app:rank}
%\begin{center} {\large \sc A. Proof of Proposition \ref{pro:rank}\\} \end{center}
%
%\begin{proof}
 The matrix $M$ to be reconstructed is factorized as Eq.~(\ref{eq:OriginalMatrix}), 
where $\Sigma=\diag(\Sigma_1,\ldots,\Sigma_r)$ is a diagonal matrix of the singular values.
 We start from following key lemma.
\begin{lemma}\label{lem:singularvalues}
There exists a numerical constant $C>0$ such that, with high probability
\begin{eqnarray}
\left|\frac{\sigma_q}{\eps}-\Sigma_q\right| \le C \Mmax \left(\frac{\alpha}{\eps}\right)^{1/2}\;, \label{eq:singularvalues}
\end{eqnarray}
where it is understood that $\Sigma_q=0$ for $q>r$, and $\Mmax=\max\{M_{ij}\}$.
\end{lemma}
The proof of this lemma is provided in \cite{KOM09}.
Applying this result to the cost function $R(i)=\frac{\sigma_{i+1}+\sigma_1\sqrt{i/\eps}}{\sigma_i}$, 
we get the following bounds :
\begin{eqnarray*}
R(r) &\le& \frac{C \Mmax\sqrt{\alpha\eps}+(\Sigma_1+C\Mmax\sqrt{\alpha/\eps}   )\sqrt{r\eps}}{\eps\Sigma_r-C \Mmax\sqrt{\alpha\eps}} \;, \label{eq:Rbound1}\\
R(i) &\ge& \frac{\eps\Sigma_{i+1} - C \Mmax\sqrt{\alpha\eps}}{\eps\Sigma_i+C \Mmax\sqrt{\alpha\eps}} \;\;\;,\; \forall i<r \label{eq:Rbound2}\\
R(i) &\ge& \frac{(\Sigma_1-C\Mmax\sqrt{\alpha/\eps})\sqrt{r\eps}}{C \Mmax\sqrt{\alpha\eps}} \;\;\;,\; \forall i>r \label{eq:Rbound3}
\end{eqnarray*}
Let, $\beta=\Sigma_1/\Sigma_r$ and $\xi=(\Mmax\sqrt{\alpha})/(\Sigma_1\sqrt{r})$. 
After some calculus, we establish that for 
\begin{eqnarray*}
     \eps > C\,r\,\max\left\{\;\xi^2\;,\;\xi^4\beta^2\;,\;{\beta^4}\;\right\} \;,
\end{eqnarray*}
we have the desired inequality $R(r)<R(i)$  for all $i\neq r$. 
This proves the remark.

%
%============================================================
%
\section*{Acknowledgment}

We thank Andrea Montanari for stimulating discussions and 
helpful comments on the subject of this paper.
This work was partially supported by
a Terman fellowship, the NSF CAREER award CCF-0743978
and the NSF grant DMS-0806211.

%
%*****************************************************************
%
%\bibliographystyle{amsalpha}
\bibliographystyle{plain}

\bibliography{MatrixCompletion}

% that's all folks
\end{document}